\documentclass[12pt]{article}
\usepackage{amsmath, amsthm, amssymb}
\usepackage{amsfonts}
\usepackage{epsfig}
\usepackage{graphics}
\usepackage{subfigure}
\usepackage{graphicx}
\usepackage{color}
\usepackage{epstopdf} 
\usepackage{float} 
\usepackage{mathrsfs} 
\usepackage{relsize} 

\counterwithin*{equation}{section}
\newcommand{\proofend}{\qed}
\newtheorem{theorem}{Theorem}
\newtheorem{corollary}{Corollary}
\newtheorem{lemma}{Lemma}
\newtheorem{remark}{Remark}
\newtheorem{example}{Example}
\newcommand{\bq}{\begin{equation}}
	\newcommand{\eq}{\end{equation}}
\newcommand{\bqr}{\begin{eqnarray}}
	\newcommand{\eqr}{\end{eqnarray}}
\newcommand{\bqrn}{\begin{eqnarray*}}
	\newcommand{\eqrn}{\end{eqnarray*}}

\begin{document}
\begin{center}
{\large \bf Analysis of solutions of some multi-term\\ fractional Bessel equations\\[3mm]} 
	{\normalsize Pavel B. Dubovski, Jeffrey A. Slepoi\footnote{Stevens Institute of Technology\\
			Hoboken, New Jersey, USA\\
			pdubovsk@stevens.edu;\
			jslepoi@stevens.edu}\\
	\small Published in {\it Fractional Calculus and Applied Analysis}, volume {\bf 24},
	No. 5 (2021), pp. 1380–1408,\\ 
	DOI: 10.1515/fca-2021-0059 
}
\end{center}

\noindent {\small {\bf Abstract}\\
We construct the existence theory for generalized fractional Bessel differential equations and find
the solutions in the form of fractional or logarithmic fractional power series. We figure out the cases when the series solution is unique, non-unique, or does not exist. The uniqueness theorem in space $C^p$ is proved for the corresponding initial value problem. We are concerned with the following homogeneous generalized fractional Bessel equation
\vskip -7pt
\begin{equation*}
	\sum_{i=1}^{m}d_i\, x^{\alpha_i}D^{\alpha_i} u(x) + (x^\beta - \nu^2)u(x)=0, \ \  \alpha_i > 0, \beta > 0,
\end{equation*}
which includes the standard fractional and classical Bessel equations as particular cases. Mostly, we consider fractional derivatives in Caputo sense and construct the theory for positive coefficients $d_i$.
Our theory leads to a threshold admissible value for $\nu^2$, which perfectly fits to the known results. Our findings are supported by several numerical examples and counterexamples that justify the necessity of the imposed conditions.
The key point in the investigation is forming proper fractional power series leading to an algebraic characteristic equation. Depending on its roots and their multiplicity/complexity, we find the system of linearly independent solutions.

\medskip


\smallskip

}

\section{Introduction}\label{sec:1}

\setcounter{section}{1}

The Bessel functions of the first $J_{\nu}(x)$ and second kind $Y_{\nu}(x)$ are linearly independent solutions of the classical Bessel equation of mathematical physics
\vskip -12pt
\begin{equation}\label{BesselEq}
x^2 u'' + x u' + (x^2 - \nu^2)u =0.
\end{equation}
\vskip -2pt \noindent
A generalization of the 2nd order operator in this equation, that can be written also in the form
$x^2 B_{\nu} = x^2 \left(x^{\nu-1} \dfrac{d}{dx}  x^{-2\nu+1} \dfrac{d}{dx} x^{\nu}\right)$, 
 leads to the hyper-Bessel operators of arbitrary integer order:
\begin{eqnarray} \label{hB}
&& B = x^{\alpha_0}\frac{d}{dx}x^{\alpha_1}\frac{d}{dx}\ldots x^{\alpha_{m-1}}\frac{d}{dx} x^{\alpha_m}\\
&& = x^{-\beta}\left(x^m \frac{d^m}{dx^m}+a_1x^{m-1} \frac{d^{m-1}}{dx^{m-1}}+\ldots+a_{m-1}xD+a_m\right),\ \  x>0, \nonumber
\end{eqnarray}
\vskip -1pt \noindent 
where $m-(\alpha_0+\alpha_1+\ldots \alpha_m)=\beta>0$.
These operators were introduced 
by Dimovski \cite{Di1966} in 1966 and studied in a series of his papers developing operational calculi for their linear right inverse integral operators. Next, Dimovski and Kiryakova 
 continued the studies on operators \eqref{hB}.
Especially, Kiryakova devoted all Chapter 3 of the monograph  \cite{Kir1994} to detailed theory of the hyper-Bessel operators and to solutions of hyper-Bessel differential equations in explicit form by means of the hyper-Bessel functions.
In fact, the operators \eqref{hB}, the coresponding hyper-Bessel integral operators and the representations of their fractional powers in terms of integral operators involving Meijer's $G$-functions, gave the hint for the theory of the generalized fractional calculus, developed in \cite{Kir1994}.
A brief review of the theory of generalized fractional calculus as arising on the base of the hyper-Bessel operators can be found in \cite{Kir2014}.

\vskip 2pt

The most popular definitions for differentiation of fractional (non-integer) order $\alpha>0$
are the Riemann-Liouville derivative, 
defined as 
\vskip -12pt
\[
D_R^\alpha u(x)= \frac{d^n}{dx^n}\left( \frac1{\Gamma(n-\alpha)}\int_0^x (x-t)^{n-\alpha-1} u(t)dt \right)
= \frac{d^n}{dx^n}\, R^{n-\alpha} u(x),
\]
\vskip -1pt \noindent 
where $n=\lceil{\alpha}\rceil$, and the Caputo derivative,
\vskip -10pt
\[
D^\alpha u(x)= \frac1{\Gamma(n-\alpha)}\int_0^x (x-t)^{n-\alpha-1} u^{(n)}(t)dt,
\]
both based on the Riemann-Liouville fractional integral  $R^{n-\alpha}$ but with interchanged order of integer order differentiation $D^n$.
Extensive information on the theory of fractional calculus can be found in books like \cite{Samko, Podlubny, Kilbas} and recently, in many others.

Let us briefly comment some recent works on the close topic for fractional order generalizations of the Bessel equation.

In \cite{Okrasinski} 
Okrasi$\acute{\text{n}}$ski and P\l{}ociniczak  
considered a modified Bessel equation containing Riemann-Liouville fractional derivatives
\vskip -14pt
\[
x^\alpha D^\alpha \left(x^\beta D^\beta u(x)\right)-\left(x^{2\mu} + \nu^{2\mu}\right)u(x) =0,
\ \   0<\alpha,\beta,\mu \le 1.
\]
\vskip -4pt \noindent
The operator $x^\alpha D^\alpha x^\beta D^\beta$ in 
this equation is close to the 
hyper-Bessel operators \eqref{hB}.
The authors searched for solutions in the form of fractional series and estimated their asymptotic behavior.

A natural extension of the classical Bessel equation (\ref{BesselEq}) in terms of Caputo fractional derivatives
\vskip -10pt
\begin{equation}\label{BesselFrDE}
x^{2\alpha}D^{2\alpha}u(x)+x^\alpha D^\alpha u(x)+(x^{2\alpha}-\nu^2)u(x)=0, \ \  \alpha \in (0,1]
\end{equation}
was analyzed by Rodrigues, Viera and Yakubovich \cite{Rodrigues}, where a solution of \eqref{BesselFrDE}
in a form of series 
was identified for some specific values of $\nu$ depending on $\alpha$.


For other results on fractional order extensions of the Bessel equation \eqref{BesselEq}, when rewritten in the form
$L_{\nu} u(x) + u(x) = 0$,  we note also an operational calculus for the fractional powers of the Bessel operator $L_{\nu}$, developed by Bengochea \cite{Bengochea}. Namely, by means of these tools, the Cauchy IVP
\vskip -12pt
$$
w \left(L_{\nu}^{\alpha}\right)\, u(x) = g(x),\ \  L_{\nu}^{k \alpha} (0) = 0, \ k=0,1,...,n-1,
$$
was solved, where $w$ is a polynomial with complex coefficients of degree $n$, $\alpha \in {\mathbb C}\setminus \{0\}$.

Another work belongs to Al-Musalhi, Al-Salti and Karimov \cite{AlMusalhi}, where the authors consider direct and inverse source problems for a fractional diffusion equation with Caputo type hyper-Bessel differential operator of the form
\vskip -15pt
$$
\left(t^{\theta} {\frac {\partial} {\partial t}} \right)^{\alpha} u(x,t) - u_{xx} (x,t) = F (x,t).
$$
Naturally, the Mittag-Leffler functions appear as tools in this study.

Next, the wider class of the multi-index Mittag-Leffler functions (for details, see \cite{Kir2010})
was explored and related to fractional analogues of the hyper-Bessel differential operators \eqref{hB} of the form
$$
\mathcal{D} = x^{\alpha_0} D^{\delta_1} x^{\alpha_1} D^{\delta_2} x^{\alpha_2} \dots D^{\delta_{m-1}} x^{\alpha_{m-1}} D^{\delta_m} x^{\alpha_m},
$$
which are generalized fractional derivatives of multi-order $(\delta_1,...,\delta_m)$  in the sense of the generalized fractional calculus \cite{Kir1994}.
Namely, a multi-index Mittag-Leffler function was proven to solve a fractional analogue of the hyper-Bessel equation,
$ \mathcal{D}\, u(x) = \lambda\, u(x)$, see e.g. \cite{AliKir}.

\bigskip 

In this paper we expand further equation (\ref{BesselFrDE}) and investigate the following generalized Bessel equation
	\begin{equation}\label{GenBesselEqn}
		\sum_{i=1}^{m}d_i\,x^{\alpha_i}D^{\alpha_i} u(x) + (x^\beta - \nu^2)u(x)=0, \ \  x > 0,
	\end{equation}
where $\alpha_i, \beta$ are positive for all $i=1,...,m$, and the derivatives are used in the mentioned Caputo sense:
\vskip -10pt 
\begin{equation}\label{CaputoFrD} 
D^{\alpha}u(x)=\frac{1}{\Gamma(n-\alpha)}\int_{0}^{x} (x-t)^{n-\alpha-1}\, 
u^{(n)}(t) dt 
\end{equation}
Here $n-1 < \alpha < n,\, n \in \mathbb{N}$ and $\Gamma(n-\alpha)$ is the Gamma function.
If $\alpha_i = n_i$ in   \eqref{GenBesselEqn}, then the derivatives are 
considered to be  regular integer order derivatives. 
We prove the existence theorem, 
the uniqueness theorem for the initial value problem,
and provide several supporting numerical examples which vividly support our results. 

\medskip %

The results for equation \eqref{GenBesselEqn} are obtained for a broad class of parameters $\alpha_i$ and $\nu$.
They include as particular cases some results of \cite{Rodrigues} and \cite[Sect. 3.4]{Kir1994}.
For example, if we deal with derivatives of integer orders $\alpha_i$
in our equation (\ref{GenBesselEqn}), then we arrive at the eigenvalue problem of the form
 $B\,u(x)=\lambda\, u(x)$, which is solved in \cite{Kir1994} (see also \cite{KirMcBride}) and where the fundamental system of solutions has been found in terms of the hyper-Bessel functions of Delerue \cite{Delerue}.
It is worth mentioning that the Cauchy problem for the non-homogeneous hyper-Bessel differential equations
\vskip -10pt 
\begin{equation}
B\, u(x)=\lambda\, u(x)+f(x) \label{hyper1}
\end{equation}%
\vskip -3pt \noindent
was also resolved 
in \cite[Sect. 3.4]{Kir1994}, for a short exposition see in \cite{KirMcBride}.

\smallskip

We should point out that in this paper we investigate only the homogeneous equation of fractional Bessel type.

\medskip
In Section \ref{sectGenFrac} we construct the solution for the generalized fractional Bessel-type equation in the form of fractional power series.
The key point of our approach is to determine power parameter $\gamma$, which must work for all terms in the fractional power series regardless of the coefficients in the original differential equation. The remarkable phenomenon is that such a value of $\gamma$ exists and satisfies a sophisticated algebraic equation (characteristic equation). 
In Section \ref{sectAddlAnalZeros} we construct the existence theory in the space of continuously differentiable functions. 

In the next two sections we consider the scenario when the characteristic equation possesses repeated roots with multiplicity two or more. In these cases we also could construct solutions using series approach but the fractional power series switch to the logarithmic fractional series.

In Section \ref{sectUniq} we prove the uniqueness theorem for the corresponding initial value problem in the space $C^p$ of continuously $p$-differentiable functions. 

Section \ref{ExamplesSection} provides the numerical examples supporting the constructed theory. The examples
cover the following principal cases:

(1) The number of roots $\gamma$ of the above mentioned key algebraic equation coincides with the order
of the original equation with integer derivatives, and we become able to construct the entire fundamental system of linearly independent solutions to generalized fractional Bessel equation.

(2) The value of parameter $\nu^2$ does not satisfy the condition in the constructed existence theory, and the equation has no fractional power series solution. Therefore, the threshold for $\nu^2$ found in our existence theorem, is essential.

(3) The value of parameter $\nu^2$ meets the threshold condition and a fractional power series solution exists. 

(4-6) The assumptions of the constructed existence theory are broken and the positivity condition on the
coefficients $d_i$ fails, then the results outside our theory become unpredictable, and we obtain either zero, one, or many solutions.

(7) The existence of multiple series solutions in case when the highest integer derivative is above the derived threshold.

(8) The samples of linear independent solutions corresponding to the same multiple root are presented.  

\vskip 2pt 

Thus, the provided examples and counterexamples cover all cases that either fit or do not fit the conditions of the presented theory and support, from different viewpoints, the necessity of our results.

\vspace*{-1pt} 

\section{Construction the fractional series solutions}\label{sectGenFrac}

\setcounter{section}{2}

Let us consider equation
\vskip -11pt
\begin{equation}\label{GenFrac}
\sum_{i=1}^{m}d_i x^{\alpha_i}D^{\alpha_i}u(x)+(x^\beta-\nu^2)u(x)=0 \text {,  where } d_i \in \mathbb{R},
\ \alpha_i>0,\ \beta>0.
\end{equation}

%
If we assume that the solution can be represented as a series\\
$\displaystyle u(x)=\sum_{n=0}^{\infty}c_n x^{\gamma+\beta n}$,
then the substitution at $c_0=1$ yields
 \vskip -10pt %
\begin{eqnarray}\label{GenSerSoln}
&&\!\!\!\!\!\!\!\!\!\!\!\!\!c_0 \cdot x^\gamma\left(\sum_{i=1}^{m}\frac{d_i\Gamma(1+\gamma)}{\Gamma(1+\gamma-\alpha_i)}-\nu^2\right)+\sum_{n=0}^{\infty}c_{n+1} x^{\gamma+\beta (n+1)}  \nonumber \\
&&\!\!\!\!\!\!\!\!\!\!\!\!\!\times\left(\sum_{i=1}^{m}\frac{d_i\Gamma(1+\gamma+\beta (n+1))}{\Gamma(1+\gamma+\beta (n+1)-\alpha_i)}-\nu^2\right)
+\sum_{n=0}^{\infty}c_n x^{\gamma+\beta (n+1)}=0.
\end{eqnarray}

In order to zero the coefficient at $x^\gamma$ in equation \eqref{GenSerSoln}, parameter $\gamma$ should satisfy the following characteristic equation
\begin{equation}\label{CondForC0} 
\sum_{i=1}^{m}\frac{d_i\cdot \Gamma(1+\gamma)}{\Gamma(1+\gamma-\alpha_i)}=\nu^2.
\end{equation}
\newpage 
\noindent 
The two series in \eqref{GenSerSoln} undergo the cancellation if
\vskip -12pt 
\begin{eqnarray} \label{Eqnforc_n} %
   c_{n+1}&=&\frac{-c_n}{\displaystyle\sum_{i=1}^{m}\frac{d_i\Gamma(1+\gamma+\beta (n+1))}{\Gamma(1+\gamma+\beta (n+1)-\alpha_i)}-\nu^2} \nonumber \\
  &=&\frac{(-1)^{n+1}c_0}{\displaystyle\prod_{k=1}^{n+1}\left(\displaystyle\sum_{i=1}^{m}\frac{d_i\Gamma(1+\gamma+\beta k)}{\Gamma(1+\gamma+\beta k-\alpha_i)}-\nu^2\right)}, \ \  n \ge 0.
 \end{eqnarray}
Hence, the solution to equation \eqref{GenFrac} can be expressed as
\vskip -12pt
\begin{eqnarray}\label{SolnGenBessel} 
u(x)=c_0\cdot{\mathlarger{\mathlarger{\mathlarger\sum_{n=0}^{\infty}}}}\frac{(-1)^n x^{\gamma+\beta n}}
{\displaystyle\prod_{k=1}^{n}\left(\displaystyle\sum_{i=1}^{m}\frac{d_i\Gamma(1+\gamma+\beta k)}{\Gamma(1+\gamma+\beta k-\alpha_i)}-\nu^2\right)}
\end{eqnarray}
\vskip-3pt\noindent
where $c_0$ is any constant.

\begin{lemma}\label{lemma0} Let $\gamma$ be a root of the characteristic equation (\ref{CondForC0}) and there be no positive integer value of $k$ such that $\gamma+\beta k$ is another root of (\ref{CondForC0}).
Then series \eqref{SolnGenBessel} is pointwise convergent for all $x > 0$.
\end{lemma}

\proof
Using the Stirling asymptotic approximation
\begin{equation*}
\lim_{k\to\infty}\frac{\Gamma(k+\alpha)}{\Gamma(k)k^\alpha}=1,
\end{equation*}
\vskip -4pt \noindent	
we obtain
\vskip -12pt
\begin{equation*}
\frac{\Gamma(1+\gamma+\beta k)}{\Gamma(1+\gamma+\beta k-\alpha_i) } \sim (\beta k)^{\alpha_i}
\end{equation*}
for $k\to \infty$.
\noindent Then there exists a positive constant $a>0$ such that for sufficiently big values of $k$
\vskip -10pt
\begin{equation*}
\left|\sum_{i=1}^{m}{\left({\frac{d_i\Gamma(1+\gamma+\beta k)}{\Gamma(1+\gamma+\beta k-\alpha_i) }}\right)}\right|
\geq a (\beta k)^{\alpha_{\max}},
\end{equation*}
\vskip -2pt \noindent	
where $\alpha_{\max}=\max\{\alpha_i\} > 0$.
Then for sufficiently big $n$
\vskip -10pt
\[
\prod_{k=1}^{n}\left|\displaystyle\sum_{i=1}^{m}\frac{d_i\Gamma(1+\gamma+\beta k)}{\Gamma(1+\gamma+\beta k-\alpha_i)}-\nu^2\right| \geq a_1 (\beta^n n!)^{\alpha_{\max}},\ a_1>0.
\]
Consequently, series \eqref{SolnGenBessel} converges.
It is worth pointing out that we assume that $\gamma_2=\gamma+\beta k$ is not a root for equation (\ref{CondForC0}). Otherwise, the lower value of $\gamma$ is not a valid root since the above product becomes zero. More details on this issue are given in Remark \ref{remark1_}.
\proofend 

\medskip

In general, we cannot guarantee that algebraic equation \eqref{CondForC0} has a real root $\gamma$. Some combinations of orders of fractional derivatives $\alpha_i$ and coefficients $d_i$ don't generate solutions. However, for positive coefficients $d_i$ and sufficiently broad class of parameter $\nu^2$, the existence and uniqueness of the series solution can be guaranteed as we show below.

\vspace*{-3pt} 

\section{Existence and uniqueness}\label{sectAddlAnalZeros}

\setcounter{section}{3}
\setcounter{equation}{0}\setcounter{theorem}{0}

To understand whether equation \eqref{CondForC0} is solvable, we first analyze the case of all integer
 $\alpha_i \in \mathbb{N}$.

\vspace*{-3pt}

\begin{lemma}\label{lemma1} 
	There exists a solution for the auxiliary equation \eqref{CondForC0} with $\alpha_i \in \mathbb{N} \text{ and } d_i > 0$, $i=1,...,m$.
\end{lemma}

\vspace*{-2pt}%

\proof
Let us consider the left-hand side of equation \eqref{CondForC0}
\vskip -10pt
\begin{equation}\label{FunctionG}
G(\gamma)=\sum_{i=1}^{m}\frac{d_i\cdot \Gamma(1+\gamma)}{\Gamma(1+\gamma-\alpha_i)}.
\end{equation}
\vskip -2pt \noindent %
In this case, function $G(\gamma)$ is a continuous function for $\gamma > -1$ because the numerator $\Gamma(1+\gamma) $ is always positive, and the denominator $\Gamma(1+\gamma-\alpha_i)$ cannot be zero. It can become infinite when $1+\gamma-\alpha_i$ take nonpositive integer values. Particularly, it happens at $\gamma=0$ when $G(0)=0$.

In addition, $G(\gamma)$ infinitely increases to infinity as $\gamma \to \infty$. Really, the application of Stirling formula yields
\vskip -10pt
\[
\displaystyle \frac{\Gamma(1+\gamma)}{\Gamma(1+\gamma-\alpha_i)} \sim \gamma^{\alpha_i}.
\]
\vskip -2pt \noindent
Consequently, since half-line $[0,+\infty)$ is included in the range of continuous function $G$, then
there exists at least one root of equation $G(\gamma)=\nu^2$.
\proofend 

\begin{remark}\label{remark1}
	Existence of the solution for auxiliary equation \eqref{CondForC0} does not guarantee the existence of a series solution because	Caputo derivative may not exist for a particular $\gamma$. We address this question in Lemmas \ref{lemma2}, \ref{lemma3} and Theorem \ref{Thrm1}.
	\end{remark}

\vspace*{-16pt}

\begin{remark}\label{remark1_}
	If characteristic equation (\ref{CondForC0}) has two or more roots distant by $\beta n$ with positive integer $n$, then one of the coefficients $c_n=\infty$ and $\gamma_i-\gamma_j = \beta n$. That implies that the smaller root $\gamma_j$ must not be taken into account, it does not represent a valid solution for equation \eqref{GenFrac}. That's why this condition is included in the statement of Lemma \ref{lemma0}. This issue is clarified by the example below.
	\end{remark}

\vspace*{-16pt}	%

	\begin{example} 
For differential equation $x^4 u^{(4)}+xu=0$,
	which corresponds $\nu=0, \beta = 1$, $m=1$, the auxiliary equation \eqref{CondForC0} becomes
\vskip-8pt
	\begin{equation*}
	\frac{\Gamma(1+\gamma)}{\Gamma(1+\gamma-4)}=0.
	\end{equation*}
	This equation has four roots: $\gamma_1=0, \gamma_2=1, \gamma_3=2, \gamma_4=3$.  For $\gamma_1, \gamma_2, \gamma_3$ denominator in $c_1$ is
	$\displaystyle\frac{\Gamma(1+\gamma+1)}{\Gamma(1+\gamma+1-4)}=0$, which makes $c_1=\infty$, and the series solution does not exist.
	Only for $\gamma_4=3$ we obtain valid result
\vskip -10pt
	$$ 
c_1=
	-\displaystyle\frac{c_0}{\Gamma(5)}=-\frac{c_0}{4!}; \, c_2=\displaystyle\frac{c_0\Gamma(2)}{\Gamma(5)\Gamma(6)}=\frac{c_0}{4!5!}; \, c_3=-\frac{c_0\cdot2!}{4!5!6!}; \, c_4=\frac{c_0\cdot2!3!}{4!5!6!7!}. $$ 
\vskip -2pt \noindent
	In general,
\vskip -12pt
	\[
	c_n=\frac{c_0(-1)^n 2! 3!}{n!(n+1)! (n+2)! (n+3)!}, \ \, n \ge 4.
	\]
	Thus, up to a constant multiplier,
	\[
	u(x)=x^3-\frac{x^4}{4!}+\frac{x^5}{4!5!}-\frac{x^6\cdot 2!}{4!5!6!}+\sum_{n=4}^\infty \frac{(-1)^n 2! 3!}{n!(n+1)! (n+2)! (n+3)!} x^{n+3}
	\]
	is the only solution in the series form.
	\end{example}
\vspace*{-2pt}%

{\sc Notation.}
Let $n_{\max}$ be the integer ceiling for the highest non-integer derivative,
i.e., $n_{\max}=\max{\{n_i\}}$ and $n_{\min}$ be the integer ceiling for the lowest non-integer derivative,
i.e., $n_{\min}=\min{\{n_i\}}$ where $\alpha_i\in(n_i-1,n_i)$, $n_i \in \mathbb{N}, i=1,...,m_0$.
Here $m_0 \le m$ is the number of genuinely fractional derivatives and $m-m_0$ is the number of integer derivatives in equation \eqref{GenFrac}. Let also $\displaystyle\alpha_{\max}=\max_{1\leq i\leq m}\{\alpha_i\}$ and
$p=\lceil \alpha_{\max} \rceil$.

Since $n_{\max}$ is defined by truly {\it fractional}
derivatives, then $n_{\max} \leq \lceil \alpha_{\max} \rceil$.
 In Lemma \ref{lemma2} below we consider equations, which have at least one fractional non-integer derivative.
\vspace*{-3pt}

\begin{lemma}\label{lemma2}
For the existence of the series solution for equation \eqref{GenFrac},
it is necessary that in equation \eqref{CondForC0},  $\gamma > n_{\max}-1 \ge 0$.
\end{lemma}

\vspace*{-4pt}

\proof
	We assume that at least one $\alpha_i$ is genuinely fractional,	i.e. $\alpha_i \in \mathbb{R_+}\setminus\mathbb{N}$.
If $\gamma \le n_{\max}-1$, then the Caputo derivative of order \\
$\displaystyle \alpha_{\max}=\max_{1 \le i \le m_0}\{\alpha_i\}$ of function $x^\gamma$ is divergent.
This proves Lemma \ref{lemma2}.
\proofend 

\begin{lemma}\label{lemma3}
	\begin{equation}
	g(x)=\frac{\Gamma(x)}{\Gamma(x-\alpha)}, \ \  x > \max \{0, \alpha\},
	\end{equation}
\vskip -1pt \noindent
	is a monotonically increasing function.
\end{lemma}

\vspace*{-2pt}

\proof
Let us inspect the derivative of $g(x)$:
\[
g'(x)=\frac{\Gamma(x)}{\Gamma(x-\alpha)}[\psi(x)-\psi(x-\alpha)],
\]
where $\psi$ is the digamma function, the logarithmic derivative of $\Gamma(x)$. Since the $\Gamma$-function is positive for positive arguments and $\psi$ is an increasing function on $(0,\infty)$, then $g'(x) > 0$ and $g(x)$ is a monotonically increasing function.
\proofend
\vspace*{-3pt} 

\begin{theorem}{}\label{Thrm1}
Let for all $i$, $1\leq i\leq m$, the coefficients $d_i$ be positive. Let $\nu$ satisfy the inequality
\vskip -12pt 	%
	\begin{equation}\label{TrueNuCondinTheorem}
	\nu^2\geq \nu^2_{\min} = \Gamma(n_{\max})\sum_{i=1}^{m}\frac{d_i}{\Gamma(n_{\max}-\alpha_i)}.
	\end{equation}
\vskip -4pt \noindent	
	Then in any domain $x \in [0,b]$, $0<b<\infty$, equation \eqref{GenFrac}
	has at least one solution $u\in C^p[0,b]$, which can be represented as fractional power series
	\eqref{SolnGenBessel}.
\end{theorem}

\proof
If all derivatives in \eqref{GenFrac} are integer, then $\nu$ can be any number for the solution to exist per Lemma \ref{lemma1}.  Now, let's consider the equation where at least one derivative is fractional.

	In equation \eqref{CondForC0} we want to find $\gamma$ based on $\nu$.
	In this theorem we can assume the smallest acceptable $\gamma$ and identify corresponding $\nu$.  We know from Lemma \ref{lemma2} that for the existence of the series solution we should consider $\gamma > n_{\max}-1$.  For
	$\gamma>n_{\max}-1$, let $\gamma=n_{\max}-1+\eta,  \eta > 0$. We introduce the function
	\vskip -10pt%
	\begin{equation*}
g_i(\eta)=\frac{\Gamma(1+\gamma)}{\Gamma(1+\gamma-\alpha_i)}=\frac{\Gamma(1+n_{\max}-1+\eta)}{\Gamma(1+n_{\max}-1+\eta-\alpha_i)}=
	\frac{\Gamma(n_{\max}+\eta)}{\Gamma(n_{\max}+\eta-\alpha_i)}.
	\end{equation*}
\vskip -2pt \noindent		%
	If $\alpha_i \ge n_{\max}$ then $\alpha_i \in \mathbb{N}$ and $\displaystyle\frac{d_i}{\Gamma(n_{\max}-\alpha_i)}=0$. Hence,  the terms with integer derivatives equal or above $n_{\max}$  do not contribute to the sum in \eqref{TrueNuCondinTheorem}. Then let's consider only $\alpha_i < n_{\max}$.
	In this case, based on Lemma \ref{lemma3}, $g_i(\eta)$ is a positive monotonically increasing function.
	Therefore, its minimum is attained at $\eta=0$ and is equal to
	$\displaystyle g_i(0)= \frac{\Gamma(n_{\max})}{\Gamma(n_{\max}-\alpha_i)}$.
Then the left-hand side of \eqref{CondForC0} is monotonically increasing for $\gamma > n_{\max}-1$ from
$\displaystyle \frac{\Gamma(n_{\max})}{\Gamma(n_{\max}-\alpha_i)}$ to infinity. Consequently, if the
threshold inequality \eqref{TrueNuCondinTheorem} holds, then the key algebraic equation
\eqref{CondForC0} has at least one root $\gamma$, and the fractional power series solution to \eqref{GenFrac}
exists. This solution has $p$ times continuous derivative because series \eqref{SolnGenBessel} converges,
along with its derivatives, uniformly on any compact interval $[0,b]$ to the continuous function $u\in C^p[0,b]$.
\proofend

\vspace*{-2pt}
\begin{corollary}
	In the case of standard fractional Bessel equation \eqref{BesselFrDE}: $\displaystyle\nu^2_{\min}=\frac{1}{\Gamma(1-2\alpha)}+\frac{1}{\Gamma(1-\alpha)}$.  This expression matches the value of $\nu^2$ identified in \cite{Rodrigues} for which solution to equation \eqref{BesselFrDE} was found.
	In our case it is the minimum acceptable value of $\nu^2$ for which the series solution is identified.\\
\end{corollary}

\vspace*{-16pt}

\begin{corollary}
	Classical Bessel equation \eqref{BesselEq} yields $\nu^2_{\min}=0$.
\end{corollary}	

\vspace*{-5pt}

\proof
	Formally, in the limit sense
	\begin{eqnarray}
	\nu^2_{\min}&=&\frac{\Gamma(0)}{\Gamma(-1)}+\frac{\Gamma(0)}{\Gamma(-2)}=\lim_{\epsilon\to 0}\frac{\Gamma(\epsilon)}{\Gamma(-1+\epsilon)}+\frac{\Gamma(\epsilon)}{\Gamma(-2+\epsilon)} \nonumber \\
	&=& \lim_{\epsilon\to 0}\epsilon \ln \epsilon + \lim_{\epsilon\to 0}\epsilon^2 \ln \epsilon = 0, \nonumber
	\end{eqnarray}
\vskip -2pt \noindent
	relying on the definition of the gamma function in the integral form.
\proofend	

\begin{remark}\label{RemarkOnNu}
	The substitution in \eqref{GenSerSoln} works only if $\gamma$ is large enough for Caputo derivative to exist.
	If $\gamma=0,...,n_{\min}-1, n_{\min}$ all derivatives of the first element of the series are zero and the series solution can be obtained if $\nu=0$ and $\beta \ge n_{\max}-1$ (to ensure the existence of the Caputo derivative). This is a solution that exists only if $\nu=0$. We do not analyze this simple special case in this work.
\end{remark}

\vspace*{-15pt}

\begin{theorem}{}{\rm (First uniqueness theorem)\ } \label{UniqThrm1} 
\\
	Let the conditions of Theorem \ref{Thrm1} hold. Then:
	\begin{enumerate}
		\item if $\nu = 0$ and $\beta \ge n_{\max}-1$ then at least one solution in the form of series can always be found;
		\item if $\alpha_{\max}$, the highest value of $\alpha_i$, is genuinely  fractional and $\nu_{\min} > 0$, then the series solution is unique;
		\item if $\alpha_{\max}$ is integer and $\nu_{\min} > 0$, then the series solution is unique provided that the difference between integer $\alpha_{\max}$ and the highest {\it fractional} derivative is less than {\rm 2}, i.e.,
		$n_{\max} \ge \alpha_{\max}-1$;
		\item if the difference between integer $\alpha_{\max}$ and the highest {\it fractional} derivative is more than {\rm 2}, i.e.,	
$n_{\max} < \alpha_{\max}-1$, then there may be several solutions in the form of series.
	\end{enumerate}
\end{theorem}

\proof
Statement 1 follows from Remark \ref{RemarkOnNu}.  

If $\alpha_{\max}$ is genuinely fractional and condition \eqref{TrueNuCondinTheorem} is satisfied, then \\
$\gamma > n_{\max}-1$ and each term in the function $G(\gamma)$ is monotonically increasing to $\infty$ per Lemma \ref{lemma3}. Hence, the series solution is unique which proves Statement 2. 

If for integer $\alpha_{\max}$ condition \eqref{TrueNuCondinTheorem} holds
and  $n_{\max} \ge \alpha_{\max}-1$, then all terms in the function $G(\gamma)$ are monotonically
increasing and the series solution is still unique.  This proves Statement 3. 

Finally, in the case $n_{\max} < \alpha_{\max}-1$ and valid condition \eqref{TrueNuCondinTheorem}, the solutions exist but the uniqueness is not guaranteed because the key algebraic equation \eqref{CondForC0} may have several roots ((see the non-uniqueness Example \ref{ex7}). Sometimes the solution may remain unique, that
happens, e.g., if the conditions of Remark \ref{remark1_}
hold and, consequently, the extra values of $\gamma$ are dummy and do not generate additional solutions. However, at least one solution still exists.
\proofend

\smallskip

Thus, Theorem \ref{UniqThrm1} guarantees that if condition \eqref{TrueNuCondinTheorem} for $\nu$ is satisfied,  all $d_i>0$, and either $\alpha_{\max}$ is fractional or $n_{\max} \ge \alpha_{\max}-1$ in case of
integer $\alpha_{\max}$, then there exists a unique $\gamma > n_{\max}-1$, which represents a unique series solution of equation \eqref{CondForC0}.

If at least one coefficient $d_i$ is negative, then function $G(\gamma)$  may have no zeros, one zero or multiple zeros, and the solution may not exist or be non-unique (see Examples \ref{ex4}--\ref{ex6} in Section \ref{ExamplesSection}).

\vspace*{-8pt} 

\section{Multiple root with multiplicity two}\label{SectBesselMult2}
\setcounter{section}{4}
\setcounter{equation}{0}\setcounter{theorem}{0}

Let us assume that a root $\gamma$ of the characteristic equation \eqref{CondForC0} has multiplicity two.  We know how to find one solution \eqref{SolnGenBessel}.  We need identify the second one. In this case, we search for the second solution to \eqref{GenFrac} in the following form
\vskip -12pt%
\begin{equation}\label{SolnMult2} 
	u(x)=C_2\left(\ln x\sum_{n=0}^{\infty}c_n^1 x^{\gamma+\beta n}+\sum_{n=0}^{\infty}c_n^2 x^{\gamma+\beta n}\right),
\end{equation}
\vskip -3pt \noindent
where $C_2$ is an arbitrary constant.

The first element of the row is $c_0^1 x^\gamma \ln x  + c_0^2 x^\gamma$.  It is known (e.g., \cite{Podlubny}, page 310) that for $\Re({\gamma})>-1$,
\vskip -12pt %
\begin{eqnarray}\label{LeibnizRulej=1} 
D^\alpha \left(x^\gamma \ln x \right) =\frac{x^{\gamma-\alpha}\Gamma(1+\gamma)}{\Gamma(1+\gamma-\alpha)}\left[\ln x + \psi(1+\gamma)-\psi(1+\gamma-\alpha)\right].
\end{eqnarray}
\vskip -3pt \noindent%
Here $\psi(x)$ is the digamma function.
At the multiple root the derivative of the function in characteristic equation \eqref{CondForC0} is equal to zero:
\vskip -12pt %
\begin{eqnarray}\label{Deriv1Eq0}
\begin{aligned}
&\left(\sum_{i=1}^{m}\frac{d_i\cdot \Gamma(1+\gamma)}{\Gamma(1+\gamma-\alpha_i)}\right)'
 \nonumber \\
&\phantom{=}=\sum_{i=1}^{m}d_i\left(\frac{\Gamma(1+\gamma)}{\Gamma(1+\gamma-\alpha_i)}\cdot [\psi(1+\gamma)-\psi(1+\gamma-\alpha_i)]\right) =0.
\end{aligned}
\end{eqnarray}
\vskip -4pt \noindent
Hence, taking into account \eqref{LeibnizRulej=1}, the Caputo derivative of the first element of the series becomes
\vskip -12pt
\begin{eqnarray}\label{DerivMult2}
\sum_{i=1}^{m}D^{\alpha_i}(c_0^1 x^\gamma \ln x + c_0^2 x^\gamma)
=\sum_{i=1}^{m}\frac{x^{\gamma-\alpha_i}\Gamma(1+\gamma)}{\Gamma(1+\gamma-\alpha_i)}[c_0^1 \ln x +c_0^2].
\end{eqnarray}
\vskip -3pt \noindent
Therefore characteristic equation \eqref{CondForC0} stays the same since it zeros out the factor at both $c_0^1$ and $c_0^2$. 

If we plug (\ref{SolnMult2}) into equation \eqref{GenFrac} and introduce simplification\\
$\Psi_1^i(\psi,n)=\psi(1+\gamma+\beta n) - \psi(1+\gamma+\beta n -\alpha_i)$, then we obtain
\vskip -12pt %
\begin{eqnarray}\label{PlugSolnMult2}
&&\hspace*{-10mm}\sum_{n=0}^{\infty}c_n^1 x^{\gamma+\beta n}\sum_{i=1}^{m}d_i \frac{\Gamma(1+\gamma+\beta n)}{\Gamma(1+\gamma+\beta n -\alpha_i)}[\ln x -\nu^2 \ln x + \Psi_1^i(\psi,n)] \nonumber\\
&+&\sum_{n=0}^{\infty}c_n^2 \left(x^{\gamma+\beta n}\sum_{i=1}^{m}d_i \frac{\Gamma(1+\gamma+\beta n)}{\Gamma(1+\gamma+\beta n -\alpha_i)}-\nu^2\right)\nonumber \\
&+&\sum_{n=0}^{\infty}c_n^1 x^{\gamma+\beta (n+1)}\ln x + \sum_{n=0}^{\infty}c_n^2 x^{\gamma+\beta (n+1)}=0.
\end{eqnarray}
\vskip -3pt \noindent %
Since characteristic equation \eqref{CondForC0} eliminates all terms with coefficients
$c_0^k, k=1,...m$, then
\vskip -12pt%
\begin{eqnarray}\label{ChkSolnMult2}
&&\hspace*{-10mm}\sum_{n=0}^{\infty}c_{n+1}^1 x^{\gamma+\beta (n+1)}\left(\sum_{i=1}^{m}d_i \frac{\Gamma(1+\gamma+\beta (n+1))}{\Gamma(1+\gamma+\beta (n+1) -\alpha_i)}-\nu^2\right)\ln x  \nonumber \\
&+&\hspace*{-2mm} \sum_{n=0}^{\infty}c_{n+1}^1 x^{\gamma+\beta (n+1)}\sum_{i=1}^{m}d_i
\frac{\Gamma(1+\gamma+\beta (n+1))}{\Gamma(1+\gamma+\beta (n+1) -\alpha_i)}\Psi_1^i(\psi,n+1) \nonumber\\
&+&\sum_{n=0}^{\infty}c_{n+1}^2 x^{\gamma+\beta (n+1)} \left(\sum_{i=1}^{m}d_i \frac{\Gamma(1+\gamma+\beta (n+1))}{\Gamma(1+\gamma+\beta (n+1) -\alpha_i)}-\nu^2\right)\nonumber \\
&+&\sum_{n=0}^{\infty}c_n^1 x^{\gamma+\beta (n+1)}\ln x + \sum_{n=0}^{\infty}c_n^2 x^{\gamma+\beta (n+1)}=0.
\end{eqnarray}
Matching coefficients in front of $x^{\gamma+\beta (n+1)}\ln x$  provides the relation
\vskip -13pt %
\begin{equation}\label{EqnCn1}
c_{n+1}^1=-\frac{c_n^1}{\displaystyle\sum_{i=1}^{m}d_i \frac{\Gamma(1+\gamma+\beta (n+1))}
{\Gamma(1+\gamma+\beta (n+1) -\alpha_i)}-\nu^2}.
\end{equation}
\vskip -4pt \noindent	
Once we know $c_{n+1}^1$  we can express $c_{n+1}^2$ through the following recursion
\vskip -12pt%
\begin{eqnarray}\label{EqnCn2}
c_{n+1}^2=-\frac{c_n^2+c_{n+1}^1\displaystyle\sum_{i=1}^{m}d_i
\frac{\Gamma(1+\gamma+\beta (n+1))}{\Gamma(1+\gamma+\beta (n+1) -\alpha_i)}\Psi_1^i(\psi,n+1)}
{\displaystyle \sum_{i=1}^{m}d_i \frac{\Gamma(1+\gamma+\beta (n+1))}{\Gamma(1+\gamma+\beta (n+1) -\alpha_i)}-\nu^2},
\end{eqnarray}
\vskip -4pt \noindent	%
which gives us the required coefficients.

\vspace*{-3pt}

\begin{lemma}\label{LemmaMult2} 
	Let $\gamma$ be a root of multiplicity two for characteristic equation (\ref{CondForC0}) and there is no positive integer value of $n$ such that $\gamma+\beta n$ is another root of (\ref{CondForC0}).
	Then series \eqref{SolnMult2} is pointwise convergent for all $x > 0$.
\end{lemma}

\vspace*{-4pt}

\proof
	Let us consider two sums in series \eqref{SolnMult2} separately.

\vskip 3pt

	The first sum has coefficients $c_{n+1}^1$ defined by \eqref{EqnCn1}.  Based on Lemma \ref{lemma0} it is convergent.

\vskip 3pt

	The other coefficient $c_{n+1}^2$ is represented in \eqref{EqnCn2}.  Taking into consideration that digamma function is monotonically increasing for $z>0$ and the identity
\vskip -13pt
	\begin{equation*}
	\psi(z+1)=\psi(z)+\frac{1}{z},
	\end{equation*}
\vskip -4pt \noindent
	we get
\vskip -14pt
	\begin{eqnarray}
	\Psi_1^i(\psi,n)<\psi(1+\gamma+\beta n)-\psi(1+\gamma+\beta n -\lceil\alpha_i\rceil)
	\!=\!\sum_{k=0}^{\lceil\alpha\rceil-1}\! \frac{1}{1+\gamma+\beta n-k}\nonumber
	\end{eqnarray}
\vskip -2pt \noindent
	since $\gamma > -1, \beta >0, n\in\mathbb{N}$.
	Therefore,
	\begin{eqnarray}
	\!\!|c_{n+1}^2|\!\!\!\!\!\!\!\!&\!=\!&\!\!\!\!\!\!\left| \frac{c_n^2+c_{n+1}^1\displaystyle\sum_{i=1}^{m}d_i \frac{\Gamma(1+\gamma+\beta (n+1))}{\Gamma(1+\gamma+\beta (n+1) -\alpha_i)}\Psi_1^i(\psi,n+1)}
	{\displaystyle\sum_{i=1}^{m}d_i \frac{\Gamma(1+\gamma+\beta(n+1))}{\Gamma(1+\gamma+\beta (n+1) -\alpha_i)}-\nu^2}\right| \nonumber	
\end{eqnarray} \begin{eqnarray}
&\underset{n\to\infty}{\sim}& \!\!\!\!\!\left| \frac {c_n^2+c_{n+1}^1 \cdot \displaystyle\sum_{i=1}^{m}(\beta n)^{\alpha_i} \sum_{k=0}^{\lceil\alpha_i\rceil-1}\frac{1}{1+\gamma+\beta (n+1)-k}}{\displaystyle\sum_{i=1}^{m}(\beta n)^{\alpha_i}}\right| \nonumber  \\
	&\underset{n\to\infty}{\sim}& \!\!\!\!\!\left| \frac {c_n^2+c_{n+1}^1 \cdot \displaystyle\sum_{i=1}^{m}(\beta n)^{\alpha_i-1}}{\displaystyle\sum_{i=1}^{m}(\beta n)^{\alpha_i}}\right|
	\underset{n\to\infty}{\sim} \left| \frac {c_n^2+c_{n+1}^1 \cdot (\beta n)^{\alpha_{\max}-1}}{(\beta n)^{\alpha_{\max}}}\right| \nonumber\\
	&=&
	\!\!\!\!\!\!\left| \frac {c_n^2}{(\beta n)^{\alpha_{\max}}}+\frac{c_{n+1}^1}{\beta n}\right|
	=O\left[\frac{1}{(\beta^n n!)^{\alpha_{\max}}}+\frac{1}{\beta n (\beta^{n+1} (n+1)!)^{\alpha_{\max}}}\right] \nonumber \\
	&=&\!\!\!\!\!\!\!\!O\left[\frac{1}{(\beta^n n!)^{\alpha_{\max}}}\right], \nonumber
	\end{eqnarray}
	where $\alpha_{\max}=\max\{\alpha_i\} > 0$.
	This proves Lemma \ref{LemmaMult2}.
\proofend

\vspace*{1pt} 
\section{Logarithmic fractional series solutions for the multiplicity greater than two}\label{SectSolnMultN}

\setcounter{section}{5}
\setcounter{equation}{0}\setcounter{theorem}{0}
Let us assume that the multiplicity of a root of equation \eqref{CondForC0} is
$N_m$ and all solutions related to that multiple root $\gamma$ can be represented as a series ($j=1,2,...,N_m$)
\vskip -12pt%
\begin{equation}\label{LogSerSoln} 
	u(x)=C_j\sum_{k=0}^{j-1}\left[(\ln x)^{j-1-k}\sum_{i=0}^{\infty}c_n^{k+1} x^{\gamma+\beta n}\right],
\end{equation}
\vskip -1pt \noindent
where $C_j$ is an arbitrary constant for each solution.

We showed in Section \ref{SectBesselMult2} that (\ref{LogSerSoln}) holds for $j=2$. Let us prove it for a more general case.
The following formula \cite{Samko} for $l,j\in\mathbb{N}$ holds
\vskip -12pt
\begin{equation}\label{LeibnizRule}
D^\alpha \left(x^\gamma (\ln x)^{j-1} \right) = x^{\gamma-\alpha}\sum_{l=0}^{j-1}\binom{j-1}{l}(\ln x)^{j-1-l} \cdot
\frac{\partial ^l}{\partial \gamma^l}\left( \frac{\Gamma(1+\gamma)}{\Gamma(1+\gamma-\alpha)}\right).
\end{equation}
If multiplicity of $\gamma$ is equal to $j$, then similar to \eqref{Deriv1Eq0} in case of multiplicity two, we have
\vskip -13pt
\begin{equation}\label{DerivOfRatio}
\frac{\partial ^l}{\partial \gamma^l}\left( \frac{\Gamma(1+\gamma)}{\Gamma(1+\gamma-\alpha)}\right)=0, l=1,...,j-1.
\end{equation}
\vskip -2pt \noindent
We know that for a series solution to exist we need to be able to apply recursion and therefore the factors at $c_0^k, k=1,...,j$ should be zero.
Hence, taking into account \eqref{LeibnizRule} and \eqref{DerivOfRatio} (as we did in case of multiplicity two with \eqref{LeibnizRulej=1} and \eqref{Deriv1Eq0}) the Caputo derivative of the first element of the row becomes as follows
\vskip-13pt
\begin{eqnarray}\label{DerivMultJ}
\sum_{i=1}^{m}D^{\alpha_i}\left(\sum_{k=0}^{j-1} c_0^j x^\gamma (\ln x)^k \right)
=\sum_{i=1}^{m}\left[\frac{x^{\gamma-\alpha_i}\Gamma(1+\gamma)}{\Gamma(1+\gamma-\alpha_i)}\sum_{k=0}^{j-1}c_0^{k+1}(\ln x)^k\right].
\end{eqnarray}
\vskip -2pt \noindent
Therefore characteristic equation \eqref{CondForC0} stays the same since it zeros out the factor at all
$c_0^k, k=1,...,j$ when solution \eqref{LogSerSoln} is plugged into equation \eqref{GenFrac}. 

We need to apply fractional derivatives to all terms in solution \eqref{LogSerSoln} according to equation
\eqref{GenFrac}.  Therefore it is necessary to identify derivatives of the ratio of $\gamma$-functions in
\eqref{CondForC0} as we can see in \eqref{LeibnizRule}. Equation \eqref{Deriv1Eq0} gives us the first derivative of the ratio
\vskip - 13pt
\begin{equation*}
\left(\frac{\Gamma(1+\gamma)}{\Gamma(1+\gamma-\alpha_i)}\right)' =
\frac{\Gamma(1+\gamma)}{\Gamma(1+\gamma-\alpha_i)}\cdot [\psi(1+\gamma)-\psi(1+\gamma-\alpha_i)]
\end{equation*}
\vskip -3pt \noindent
and therefore we have the following second derivative:
\vskip -13pt
\begin{eqnarray}\label{Deriv2RatioGammaFunctions}
\begin{aligned}
&\left(\frac{\Gamma(1+\gamma)}{\Gamma(1+\gamma-\alpha_i)}\right)'' =\left(\frac{\Gamma(1+\gamma)}{\Gamma(1+\gamma-\alpha_i)}\cdot [\psi(1+\gamma)-\psi(1+\gamma-\alpha_i)]\right)'\nonumber \\
&=\frac{\Gamma(1\!+\!\gamma)}{\Gamma(1\!+\!\gamma\!-\!\alpha_i)}\hspace{-1mm}\left[\left(\psi(1+\gamma)-\psi(1+\gamma-\alpha_i)\right)^2+\psi'(1+\gamma)-\psi'(1+\gamma-\alpha_i)\right].\nonumber
\end{aligned}
\end{eqnarray}
The higher derivatives look similar.  Let us for the sake of simplicity define
\begin{eqnarray}
\Psi_1(n)&=&\Psi_1^i(\psi,n)=\psi(1+\gamma+\beta n )-\psi(1+\gamma + \beta n-\alpha_i),\nonumber \\
\Psi_2(n)&=&\Psi_2^i(\psi,\psi',n)= \left(\psi(1+\gamma+\beta n) -\psi(1+\gamma+\beta n-\alpha_i)\right)^2 \nonumber \\
&&+\psi'(1+\gamma+\beta n)-\psi'(1+\gamma+\beta n-\alpha_i),\nonumber \\
\Psi_3(n)&=&\Psi_3^i(\psi,\psi',\psi'',n)= \left(\psi(1+\gamma+\beta n) -\psi(1+\gamma+\beta n-\alpha_i)\right)^3 \nonumber \\
&&+(\psi''(1+\gamma+\beta n)-\psi''(1+\gamma+\beta n-\alpha_i)) \nonumber \\
&&+3(\psi(1+x+\beta n)\psi'(1+\gamma+\beta n) +\psi(1+\gamma+\beta n-\alpha_i)\nonumber \\
&&\times\psi'(1+\gamma+\beta n-\alpha_i)-\psi(1+\gamma+\beta n)\psi'(1+\gamma+\beta n-\alpha_i) \nonumber \\
&&-\psi(1+\gamma+\beta n-\alpha_i)\psi'(1+\gamma+\beta n)), \nonumber \\
... \nonumber \\
\Psi_{k}^i(n)&=&\Psi_{k}^i(\psi,\psi',...,\psi^{(k-1)},n).
\end{eqnarray}
Then, we can say
\vskip -10pt
\begin{eqnarray}\label{DerivRatioGammasGeneralized}
\left(\frac{\Gamma(1+\gamma+\beta n)}{\Gamma(1+\gamma+\beta n-\alpha_i)}\right)^{(k)}
&\!\!\!\!\!\!\!\!\!=\!\!\!\!\!& \frac{\Gamma(1+\gamma+\beta n)}{\Gamma(1+\gamma+\beta n-\alpha_i)}\Psi_{k}^i(\psi, \psi',...,\psi^{(k-1)},n) \nonumber \\
 &\!\!\!\!\!\!\!\!\!\!=\!\!\!\!\!&\frac{\Gamma(1+\gamma+\beta n)}{\Gamma(1+\gamma+\beta n-\alpha_i)}\Psi_{k}^i(n),
\end{eqnarray}
where $\psi^{(l)}(z)=\displaystyle\frac{d^l}{d z^l}\psi(z)=\frac{d^{l+1}}{d z^{l+1}}\ln \Gamma(z)$ is the polygamma function.

\begin{remark}
	From equation \eqref{DerivRatioGammasGeneralized} we may conclude that since \\
	$G(\gamma)=\displaystyle\frac{\Gamma(1+\gamma)}{\Gamma(1+\gamma-\alpha_i)}=0$,
	and, therefore, all derivatives of  $G(\gamma)$ are also zero and each zero of $G(\gamma)$ represents zero of infinite multiplicity. However it is not the case. Really, $G(\gamma)=0$ when $\gamma=\alpha_i-k, k=1,...,\lceil \alpha_i \rceil$.
	Then,
\vskip -12pt
	\begin{eqnarray}
	G(\gamma)'&=&\frac{\Gamma(1+\alpha_i-k)}{\Gamma(1-k)}(\psi(1+\alpha_i-k)-\psi(1-k)) \nonumber \\
	&=&-\Gamma(1+\alpha_i-k)\frac{\psi(1-k)}{\Gamma(1-k)}\nonumber \\
	(\text{see \cite{Mezo}})
	&=&\Gamma(1+\alpha_i-k)\cdot (-e^{2\gamma_m (1-k)})\prod_{n=0}^{\infty}\left(1-\frac{1-k}{x_n}\right)e^{\frac{1-k}{x_n}} \ne 0,\nonumber
	\end{eqnarray}
\vskip -4pt \noindent
	where $\gamma_m$ is Euler-Mascheroni constant, $x_n, n=1,2,...$ represents zeros of the digamma function.
	 Hence, the previous conclusion that any derivative of $G(\gamma)=0$ is wrong. Actually, it shows that if equation $G(\gamma)=0$ represents characteristic equation for equation \eqref{GenFrac} with one derivative and $\nu=0$, then the multiplicity of each root $\gamma$ cannot be more than one. For multiplicity of the root $\gamma$ to be greater than one, there should be several derivatives in equation \eqref{GenFrac}. 	
\end{remark}

\vspace*{-2pt}

Now we can plug the proposed solution \eqref{LogSerSoln} into equation \eqref{GenFrac} like we did in \eqref{ChkSolnMult2}:
\vskip -14pt
\begin{eqnarray} 
&&\!\!\!\!\!\!\!\!\!\!\sum_{n=0}^{\infty}c_{n+1}^1 x^{\gamma+\beta (n+1)}\left(\sum_{i=1}^{m}d_i \frac{\Gamma(1+\gamma+\beta (n+1))}{\Gamma(1+\gamma+\beta (n+1) -\alpha_i)}-\nu^2\right)(\ln x)^{j-1}  \nonumber \\
&+& \sum_{n=0}^{\infty}c_{n+1}^1 x^{\gamma+\beta (n+1)}\sum_{i=1}^{m}d_i \frac{\Gamma(1+\gamma+\beta (n+1))}{\Gamma(1+\gamma+\beta (n+1) -\alpha_i)} \nonumber \\
&\times& \sum_{k=0}^{j-2}\Psi_{k+1}^i(n+1) \cdot (\ln x)^{j-2-k}\nonumber \\
&+&\sum_{n=0}^{\infty}c_{n+1}^2 x^{\gamma+\beta (n+1)}\left(\sum_{i=1}^{m}d_i \frac{\Gamma(1+\gamma+\beta (n+1))}{\Gamma(1+\gamma+\beta (n+1) -\alpha_i)}-\nu^2\right)(\ln x)^{j-2}  \nonumber \\
&+& \sum_{n=0}^{\infty}c_{n+1}^2 x^{\gamma+\beta (n+1)}\sum_{i=1}^{m}d_i \frac{\Gamma(1+\gamma+\beta (n+1))}{\Gamma(1+\gamma+\beta (n+1) -\alpha_i)}  \nonumber 
\end{eqnarray}
\begin{eqnarray}  \label{ChkSolnMultJ}
&\times& \sum_{k=0}^{j-3}\Psi_{k+1}^i(n+1) \cdot (\ln x)^{j-3-k}\nonumber \\
&+& ... \nonumber \\
&+&\sum_{n=0}^{\infty}c_{n+1}^{j-1} x^{\gamma+\beta (n+1)} \left(\sum_{i=1}^{m}d_i \frac{\Gamma(1+\gamma+\beta (n+1))}{\Gamma(1+\gamma+\beta (n+1) -\alpha_i)}-\nu^2\right)\ln x \nonumber  \\
&+& \sum_{n=0}^{\infty}c_{n+1}^{j-1} x^{\gamma+\beta (n+1)}\sum_{i=1}^{m}d_i \frac{\Gamma(1+\gamma+\beta (n+1))}{\Gamma(1+\gamma+\beta (n+1) -\alpha_i)}\cdot \Psi_{1}^i(n+1))\nonumber \\
&+&\sum_{n=0}^{\infty}c_{n+1}^j x^{\gamma+\beta (n+1)} \left(\sum_{i=1}^{m}d_i \frac{\Gamma(1+\gamma+\beta (n+1))}{\Gamma(1+\gamma+\beta (n+1) -\alpha_i)}-\nu^2\right)\nonumber \\
&+&\sum_{n=0}^{\infty}c_n^1 x^{\gamma+\beta (n+1)}(\ln x)^{j-1} + \sum_{n=0}^{\infty}c_n^2 x^{\gamma+\beta (n+1)}(\ln x)^{j-2} \nonumber \\
&+& ...+ \sum_{n=0}^{\infty}c_n^j x^{\gamma+\beta (n+1)}=0.
\end{eqnarray}
 \vskip -3pt \noindent %
Similar to \eqref{EqnCn1} and \eqref{EqnCn2}, we can deduce the recursive relationship for the coefficients
$c_{n+1}^k$, which need to be calculated in order as shown in \eqref{Eqnforc_nMult}
\vskip -12pt
\begin{eqnarray} \label{Eqnforc_nMult}
c_{n+1}^1&\!\!\!\!\!=\!\!\!\!\!&-\frac{c_n^1}{\displaystyle\sum_{i=1}^{m}d_i \frac{\Gamma(1+\gamma+\beta (n+1))}{\Gamma(1+\gamma+\beta (n+1) -\alpha_i)}-\nu^2},\nonumber \\
c_{n+1}^2&\!\!\!\!\!=\!\!\!\!\!&-\frac{c_n^2+\displaystyle\sum_{i=1}^{m}d_i \frac{\Gamma(1+\gamma+\beta (n+1))}{\Gamma(1+\gamma+\beta (n+1) -\alpha_i)}\cdot c_{n+1}^1\Psi_1^i(n+1)}
{\displaystyle\sum_{i=1}^{m}d_i \frac{\Gamma(1+\gamma+\beta (n+1))}{\Gamma(1+\gamma+\beta (n+1) -\alpha_i)}-\nu^2}, \nonumber \\
c_{n+1}^3&\!\!\!\!\!=\!\!\!\!\!&-\frac{c_n^3+\displaystyle\sum_{i=1}^{m}d_i \frac{\Gamma(1+\gamma+\beta (n+1))}{\Gamma(1+\gamma+\beta (n+1) -\alpha_i)}\cdot \sum_{k=1}^{2}c_{n+1}^k\Psi_{3-k}^i(n+1)}
{\displaystyle\sum_{i=1}^{m}d_i \frac{\Gamma(1+\gamma+\beta (n+1))}{\Gamma(1+\gamma+\beta (n+1) -\alpha_i)}-\nu^2}, \nonumber \\
&& ... \nonumber \\
c_{n+1}^j&\!\!\!\!\!=\!\!\!\!\!&-\frac{c_n^j+\displaystyle\sum_{i=1}^{m}d_i \frac{\Gamma(1+\gamma+\beta (n+1))}{\Gamma(1+\gamma+\beta (n+1) -\alpha_i)}\cdot \sum_{k=1}^{j-1}c_{n+1}^{k}\Psi_{j-k}^i(n+1)}
{\displaystyle\sum_{i=1}^{m}d_i \frac{\Gamma(1+\gamma+\beta (n+1))}{\Gamma(1+\gamma+\beta (n+1) -\alpha_i)}-\nu^2}. \nonumber 
\end{eqnarray}
\vskip -1pt \noindent %
In the case of Caputo derivative, all solutions that we considered have the form
$x^\gamma \cdot \ln^j x, \gamma > n_{\max}-1$ and therefore all of its derivatives up to $n_{\max}-1$ in $x$ are equal to zero.  The difference between Riemann-Liouville
and Caputo derivative for continuous function is
\vskip -10pt
\begin{equation*}
D_{R}^\alpha (x^\gamma \cdot f(x)) = D_C^\alpha \left(x^\gamma \cdot  f(x)\right) + \sum_{k=0}^{n_{\max}-1} \frac{\left(x^\gamma \cdot  f(x)\right)^{(k)}_{x=0}}{\Gamma(1+k-\alpha)}x^{k-\alpha}.
\end{equation*}
\vskip -4pt \noindent
Based on the Leibniz rule we have
\vskip -10pt
\begin{equation*}
\left(x^\gamma \cdot f(x)\right)^{(j)} = \sum_{n=0}^{j}\binom{j}{n}(x^\gamma)^{(j-n)} \cdot f^{(n)}(x),
\end{equation*}
\vskip -3pt \noindent
keeping in mind that $\gamma > n_{\max}-1$, makes $(x^\gamma)^{(j-n)}|_{x=0}=0, n=0,...,j$.
In our case $f(x)=(\ln x)^{k}, k=0,...,j-1$, and in view of $\gamma > j$, we obtain
\vskip -10pt
\begin{equation*}
\displaystyle\lim_{x \to 0}(x^\gamma)^{(j-n)}[(\ln x)^{k}]^{(n)}=0.
\end{equation*}
 Therefore
 \vskip - 11pt
\begin{equation*}
D_{R}^\alpha \left(x^\gamma \cdot (\ln x)^{j-k}\right) =  D_C^\alpha \left(x^\gamma \cdot  (\ln x)^{j-k}\right).
\end{equation*}
This confirms that within the area of existence of Caputo derivative, solution \eqref{LogSerSoln} works for the generalized Bessel equation when fractional differentiation is understood in both Riemann-Liouville and Caputo senses.

If $\gamma$ is a root of multiplicity $j$ and there is no $n\in \mathbb{N}$ such that \linebreak
$\gamma+\beta n$ is another root of the characteristic equation (\ref{CondForC0}), then series \eqref{LogSerSoln} is pointwise convergent for all $x > 0$. The proof is similar to the proof of Lemma \ref{LemmaMult2}, and the convergence of each $c_{n+1}^{k+1}$ is based on the  convergence of $c_{n+1}^{k}, k=1,...,j-1$.

\vspace*{-3pt} 

\section{Initial value problem}
\label{sectUniq}
\setcounter{section}{6}
\setcounter{equation}{0}\setcounter{theorem}{0}

In this section we prove uniqueness for the initial value problem for equation
\eqref{GenFrac} with positive coefficients $d_i$. The method is similar to \cite{Rodrigues}.

\vspace*{-2pt}

\begin{theorem}{}\label{Thrm2} {\rm (Second uniqueness theorem)}\\
	The initial value problem for fractional equation \eqref{GenFrac} with the domain
	$x \in [0,b]$  and initial conditions $u^{(j)}(0)=u^{(j)}_0, j=0,1,...,p-1,$ $\displaystyle p=\lceil \alpha_{\max}
	\rceil$, has a unique solution in space $C^p[0,b]$ provided that
\vskip -9pt	%
	\begin{equation}\label{NuCondinTheorem}
	\nu^2 > b_1^{\beta}+\sum_{i=1}^{m}q_i |d_i| b_1^{n_i}.
	\end{equation}
\vskip -5pt \noindent
	Here
\vskip - 14pt
	\begin{eqnarray}
	b_1=\max\{1,b\};\
	q_i=
	\begin{cases}
	\displaystyle \frac{1}{\Gamma(n_i-\alpha_i)(n_i-\alpha_i+1)} \!\!\!\!\!\!& ,\ n_{i}-1<\alpha_i < n_i \\
	1 & ,\ \alpha_i= n_i.
	\end{cases}
	\end{eqnarray}
\end{theorem}
\proof
	$C^p[0,b]$ is the space of $p$ times continuously differentiable on $[0,b]$ functions with norm
\vskip -13pt
	\begin{equation*}
	||u||_{C^p} = \sum_{k=0}^{p}||u^{(k)}||_C = \sum_{k=0}^{p}\max_{x\in [0,b]} |u^{(k)}(x)|.
	\end{equation*}
	The operator $T:C^p \mapsto C$ is defined as follows:
\vskip -10pt
	\begin{equation*}
	(Tu)(x)=\frac{1}{\nu^2}\left[x^\beta u(x) + \sum_{i=1}^{m}d_i x^{\alpha_i}D^{\alpha_i} u(x)\right].
	\end{equation*}
\vskip -3pt \noindent
	Then we write the equation in form $u(x)=(Tu)(x)$.  Taking into account formulas 2.4.24-26  in Corollary 2.3
	from \cite{Kilbas}, we obtain:
\vskip -10pt
	\begin{equation*}
	||D^{\alpha_i} u||_C \le k_{\alpha_i} ||u||_{C^p}, \text{ where } k_{\alpha_i} = \frac{b ^{n_i-\alpha_i}}{\Gamma(n_i-\alpha_i)(n_i-\alpha_i+1)}, \text { when } \alpha_i < n_i
	\end{equation*}
\vskip -3pt \noindent
	and
\vskip -10pt
	\begin{equation*}
	||D^{\alpha_i}u||_C=||u||_{C^p}, \text { when } \alpha_i=n_i.
	\end{equation*}
\vskip -3pt \noindent
	Therefore,
\vskip -12pt
	\begin{equation*}
	||D^{\alpha_i}u||_C \le q_i b^{n_i-\alpha_i}||u||_{C^p},
	\end{equation*}
\vskip -3pt \noindent
	and we get
\vskip -14pt
	\begin{eqnarray}
	||Tu_1-Tu_2||_C &\!\!\!=& \!\!\!\frac{1}{\nu^2}||x^\beta (u_1(x)-u_2(x)) + \sum_{i=1}^{m}d_i x^{\alpha_i}D^{\alpha_i}(u_1-u_2)(x)||_C\nonumber \\
	&\le& \!\!\!\frac{1}{\nu^2}\left(b_1^\beta||u_1-u_2||_C +\sum_{i=1}^{m} |d_i| b_1^{\alpha_i} q_i b^{n_i-\alpha_i} ||u_1-u_2||_{C^p} \right) \nonumber \\
	&\le&\!\!\! \frac{1}{\nu^2}\left(b_1^\beta +\sum_{i=1}^{m} q_i |d_i| b_1^{n_i}\right)||u_1-u_2||_{C^p}.
	\end{eqnarray}
\vskip 2pt \noindent
	Having condition \eqref{NuCondinTheorem} we conclude that $T$ is a contraction.  Hence we can apply the Banach fixed point theorem to complete the proof.
\proofend

\smallskip

\begin{remark}
	If all derivatives are integer and $b=1$ then inequality \eqref{NuCondinTheorem} becomes
	$\displaystyle \nu^2 > \sum_{i=1}^{m}|d_i|+1$.  Particularly, for the initial value problem of the classic Bessel equation \eqref{BesselEq} with integer derivatives we obtain $\nu^2 > 3$ as a sufficient condition for the uniqueness of the solution.
\end{remark}

\vspace*{-2pt} 
	
\section{Examples of uniqueness, non-uniqueness, non-existence, \break and multiplicity}\label{ExamplesSection}
\setcounter{section}{7}
\setcounter{equation}{0}\setcounter{theorem}{0}

\medskip

The numerical examples support the theorems proved in the previous sections. Computations were performed
using the substitution numerical method proposed in \cite{DS2021-1}.

\begin{figure}[H]
	\begin{center}
		\includegraphics[width=6cm]{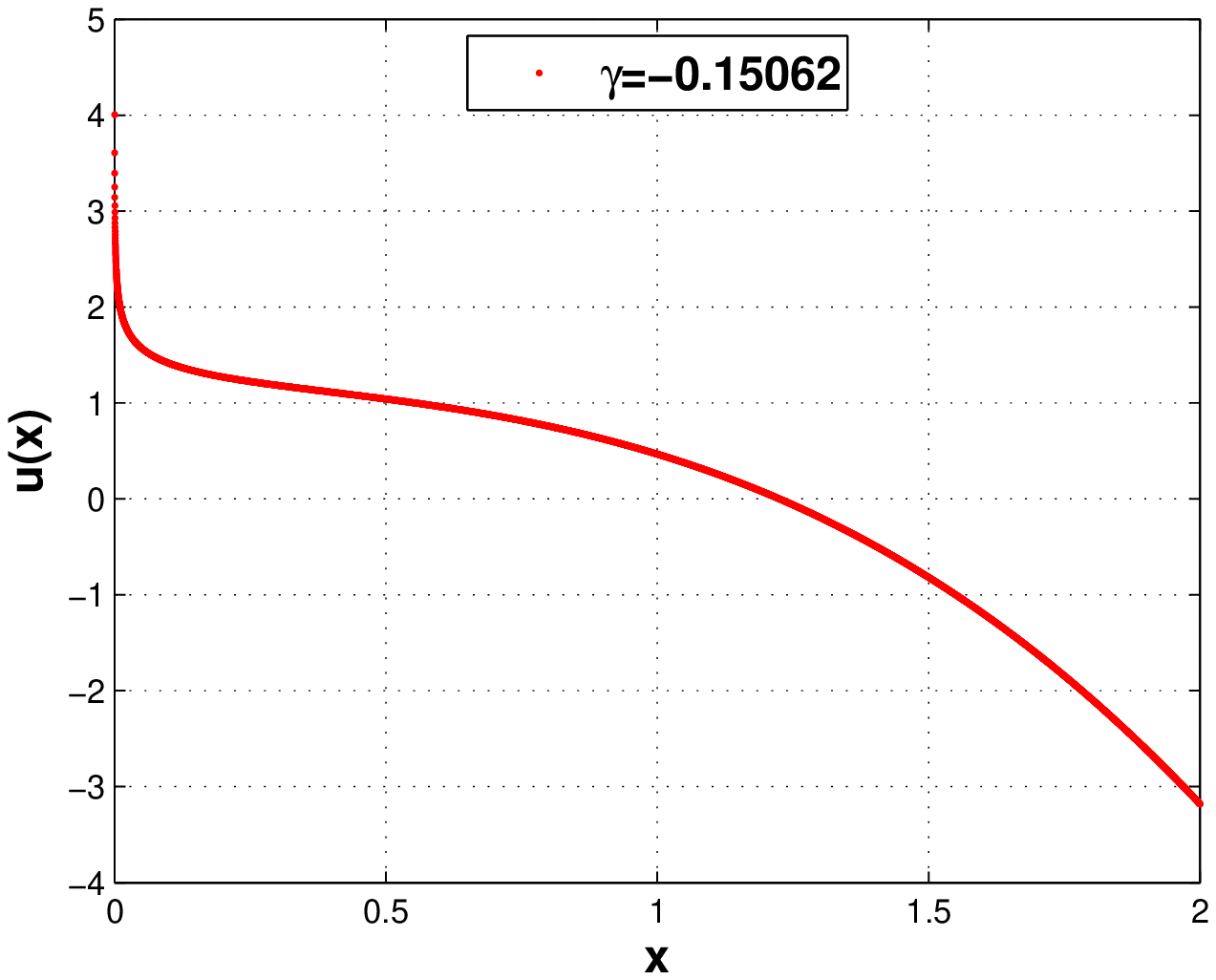}
		\includegraphics[width=6cm]{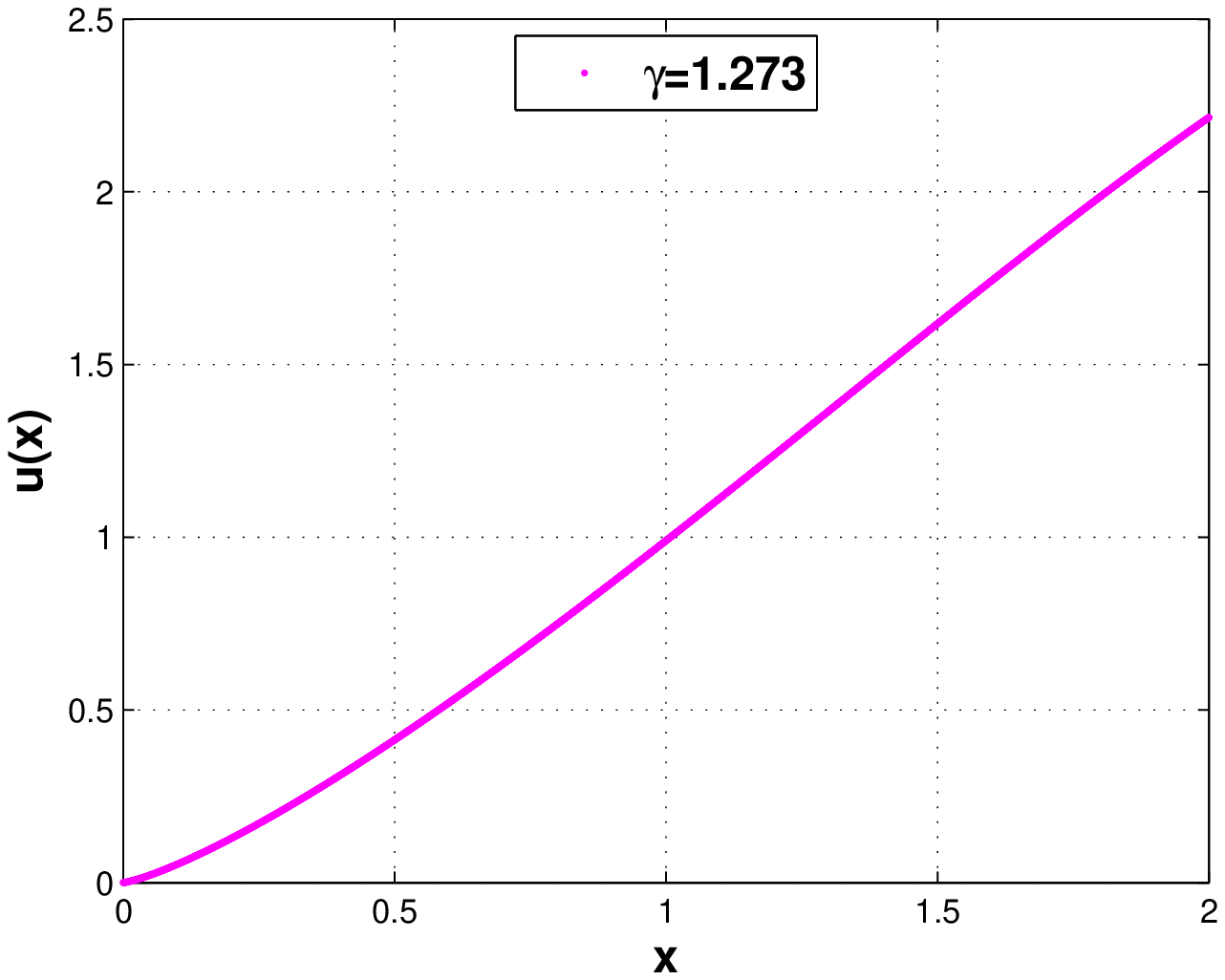}
		\includegraphics[width=6cm]{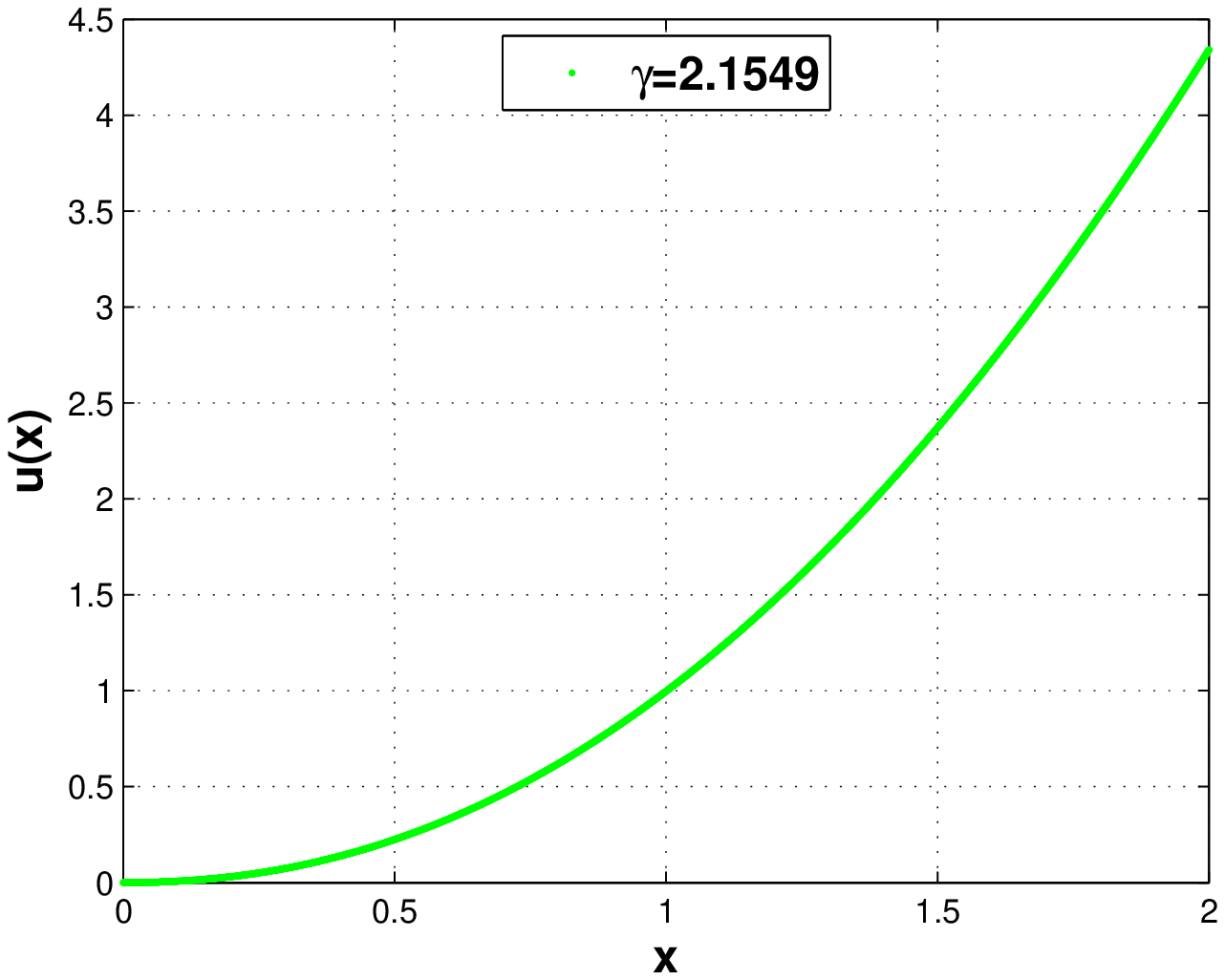}
		\includegraphics[width=6cm]{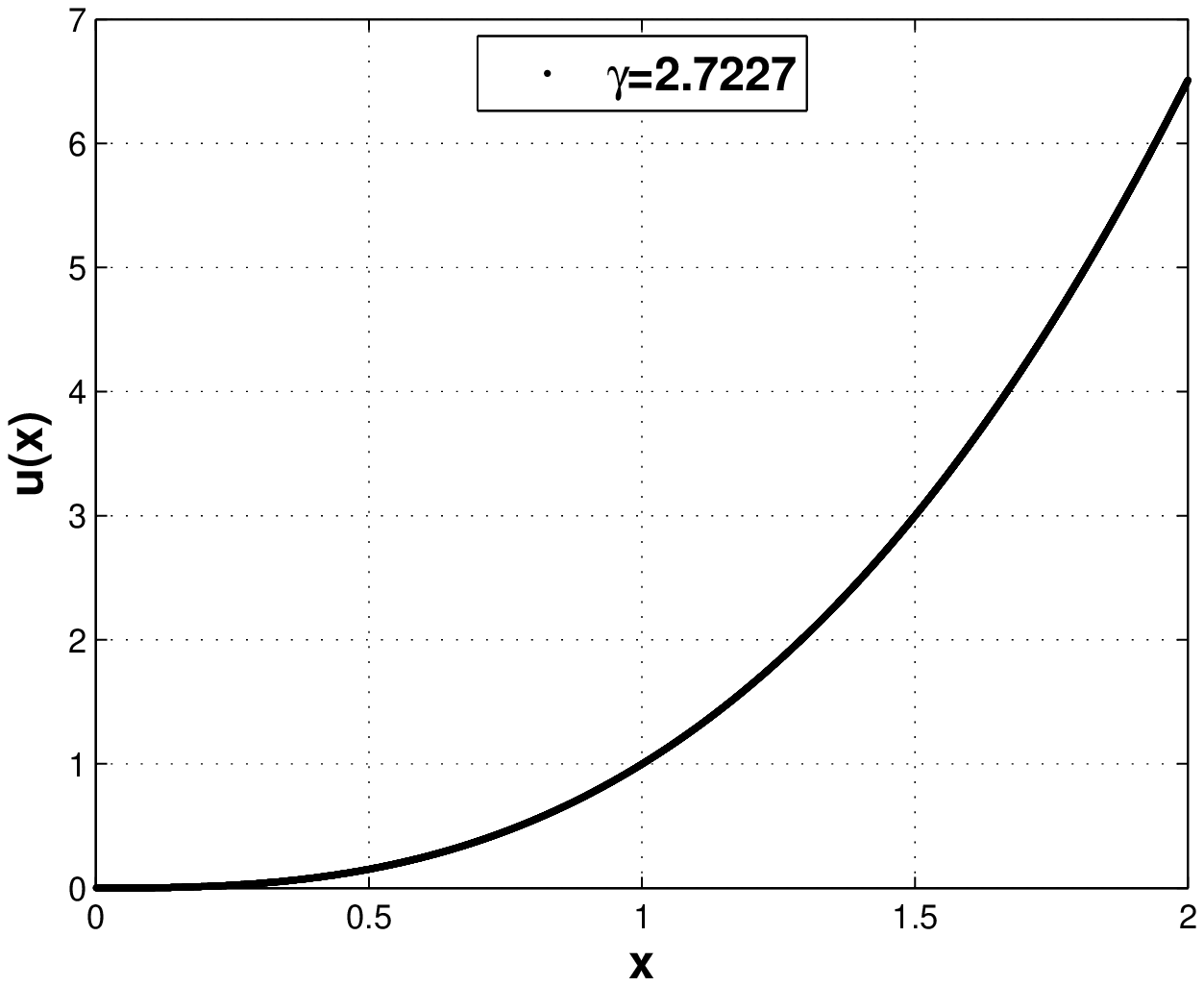}
		\caption{Solutions for equation \eqref{Example1} with integer derivatives.}	\label{SolEx1}
	\end{center}
\end{figure}

\begin{example} ($d_i>0$ and all derivatives are integer).\label{ex1}\\
Let us consider the equation
\begin{equation}\label{Example1}
2x^{4}u^{(4)}(x)+0.3x^2 u''(x) + x u'(x) + (x^{3.1}-1.5^2)u(x)=0.
\end{equation}
Here $\beta=3.1,\nu=1.5$.  Characteristic equation \eqref{CondForC0} becomes:
\begin{equation}\label{CondForC0Ex1}
\frac{2\Gamma(1+\gamma)}{\Gamma(1+\gamma-4)}-\frac{3\Gamma(1+\gamma)}{\Gamma(1+\gamma-2)}+\frac{\Gamma(1+\gamma)}{\Gamma(1+\gamma-1)}-1.5^2=0.
\end{equation}
\vskip -1pt \noindent
The graph of the expression on the left side of equation \eqref{CondForC0Ex1} is in Figure \ref{GyforEx1}. This equation has four solutions.  Please note that if $\nu$ were bigger we would have found only two solutions, but 2 series solutions exist for any $\nu$.

\vspace*{-0.7cm}
\begin{figure}[H]
	\begin{center}
		\includegraphics[width=6cm]{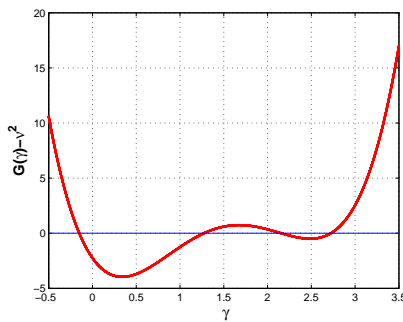}
		\caption{Function $G(\gamma)-\nu^2$ for equation \eqref{Example1}.}	\label{GyforEx1}
	\end{center}
\end{figure}

\vspace*{-0.9cm}

We use function $fzero$ in MatLab to identify zeros of the function $G(\gamma)-\nu^2$.  To be able to do it generically we need to provide the guess solution for the $fzero$ function otherwise it always finds the same zero of a function.  When the function has several zeros, the accuracy of the guess is very important.  For that, we first go through $x$-axis on the appropriate interval, in this case [-5,5], with the step of 0.001 to see where the zeros can be, i.e. $\displaystyle \left(G(\gamma_{i-1})-\nu^2\right)\cdot \left(G(\gamma_i) -\nu^2\right)\le 0$.  After that, for each interval with the zero inside we execute $fzero$ function with the guess of $\gamma_{i-1}$.  The results of the guesses in this example are  $\gamma_{g1}=-0.1505,\gamma_{g2}=1.2735$, $\gamma_{g3}=2.1545, \gamma_{g4}=2.7225$.  The actual solutions with 4 digits of precision for  equation \eqref{CondForC0Ex1} are $\gamma_1=-0.1506, \gamma_2=1.2730, \gamma_3=2.1549, \gamma_4=2.7227$.
Each of these solutions is represented in Figure \ref{SolEx1}.

\medskip 

This example shows that in some cases our method can provide the complete fundamental system of linearly independent solutions.  If the number of zeros of the function $G(\gamma)-\nu^2$ is strictly less than $n_{\max}$ then the fundamental set is incomplete.
\end{example}
\vspace*{-0.4cm}

\vspace*{-8pt} 

\begin{example} \label{ex2} (no solutions: $d_i>0$ but $\nu^2 < \nu^2_{\min}$).\\ 
Let us consider an equation with fractional derivatives
\begin{equation}\label{Example2}
x^{4.9}D^{4.9}u(x)+3.1x^{3.75}D^{3.75}u(x)+3x^{2.7}D^{2.7}u(x)+(x^{3.1}-1)u(x)=0.
\end{equation}
Here $\beta=3.1,\nu=1$.   Characteristic equation \eqref{CondForC0} becomes:
\begin{equation}\label{CondForC0Ex2}
\frac{3.1\Gamma(1+\gamma)}{\Gamma(1+\gamma-2.75)}+\frac{3\Gamma(1+\gamma)}{\Gamma(1+\gamma-2.7)}+\frac{\Gamma(1+\gamma)}{\Gamma(1+\gamma-4.9)}-1=0.
\end{equation}
There are solutions $\gamma$ to \eqref{CondForC0Ex2}: $\gamma = 0.0848, \gamma = 0.7064, \gamma = 1.9888$, none of which  can create a meaningful solution of equation \eqref{Example2} since Caputo derivative of order greater than two of $x^\gamma$  diverges.  In this example \linebreak
$n_{\max}=5$, hence, in view of Lemma \ref{lemma2} the necessary condition is $\gamma > 4$.  Based on inequality \eqref{TrueNuCondinTheorem}, for the existence of solution we need \linebreak
$\nu^2_{\min}=146.32$, but here we have $\nu^2=1$, which does not fall into the admissible interval defined in Theorem \ref{Thrm1}.
\end{example}
\vspace*{-8pt}

\begin{example} \label{ex3} (solution exists: $d_i>0, \nu^2 \ge \nu^2_{\min}$).\\ 
This example is almost verbatim Example \ref{ex2} but the value of \linebreak
$\nu^2=169 > 146.32$ falls in the admissible interval defined in Theorem \ref{Thrm1}.  Therefore it falls within the frame of the constructed theory and produces one series solution.
\begin{equation}\label{Example2a}
x^{4.9}D^{4.9}u(x)+3.1x^{3.75}D^{3.75}u(x)+3x^{2.7}D^{2.7}u(x)+(x^{3.1}-13^2)u(x)=0.
\end{equation}
Here $\gamma = 4.0934$ and the solution for equation \eqref{Example2a} can be found (see Figure \ref{Example2aGraph}).

\vspace*{-13pt} 
\begin{figure}[H]
	\begin{center}
		\includegraphics[width=6cm]{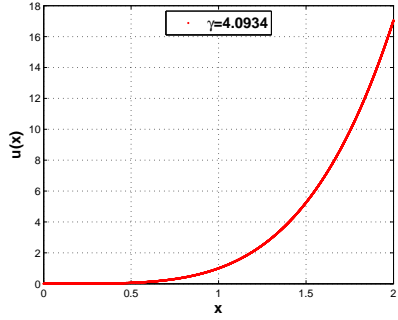}
		\caption{Solution of equation \eqref{Example2a} with $\nu=13$.}	\label{Example2aGraph}
	\end{center}
\end{figure}
\end{example}

\vspace*{-10pt}

If coefficients $d_i$ are not all positive, the above constructed theory is not applicable. We can get several solutions, no solutions or one solution.  In this case, the uniqueness and existence of a solution are not guaranteed.  The following three examples demonstrate these cases.

\vspace*{-3pt}

\begin{example} \label{ex4} (two solutions with $d_i < 0$).\\ 
Let's consider the following equation which includes negative coefficients:
\begin{equation}\label{Example2b}
-0.1x^{4.9}D^{4.9}u(x)+3.1x^{3.75}D^{3.75}u(x)-6x^{2.7}D^{2.7}u(x)+(x^{3.1}-1)u(x)=0.
\end{equation}

The calculated minimum for $\nu^2$ is negative, therefore there is no restriction on $\nu$.  In our case, with $\nu=1$, we get two acceptable values of $\gamma = 4.6583$ and  $\gamma = 22.1448$. which means that we generate two different solutions as shown in Figure \ref{Example2bGraph}.

\newpage 
\begin{figure}[H]
	\begin{center}
		\includegraphics[width=6cm]{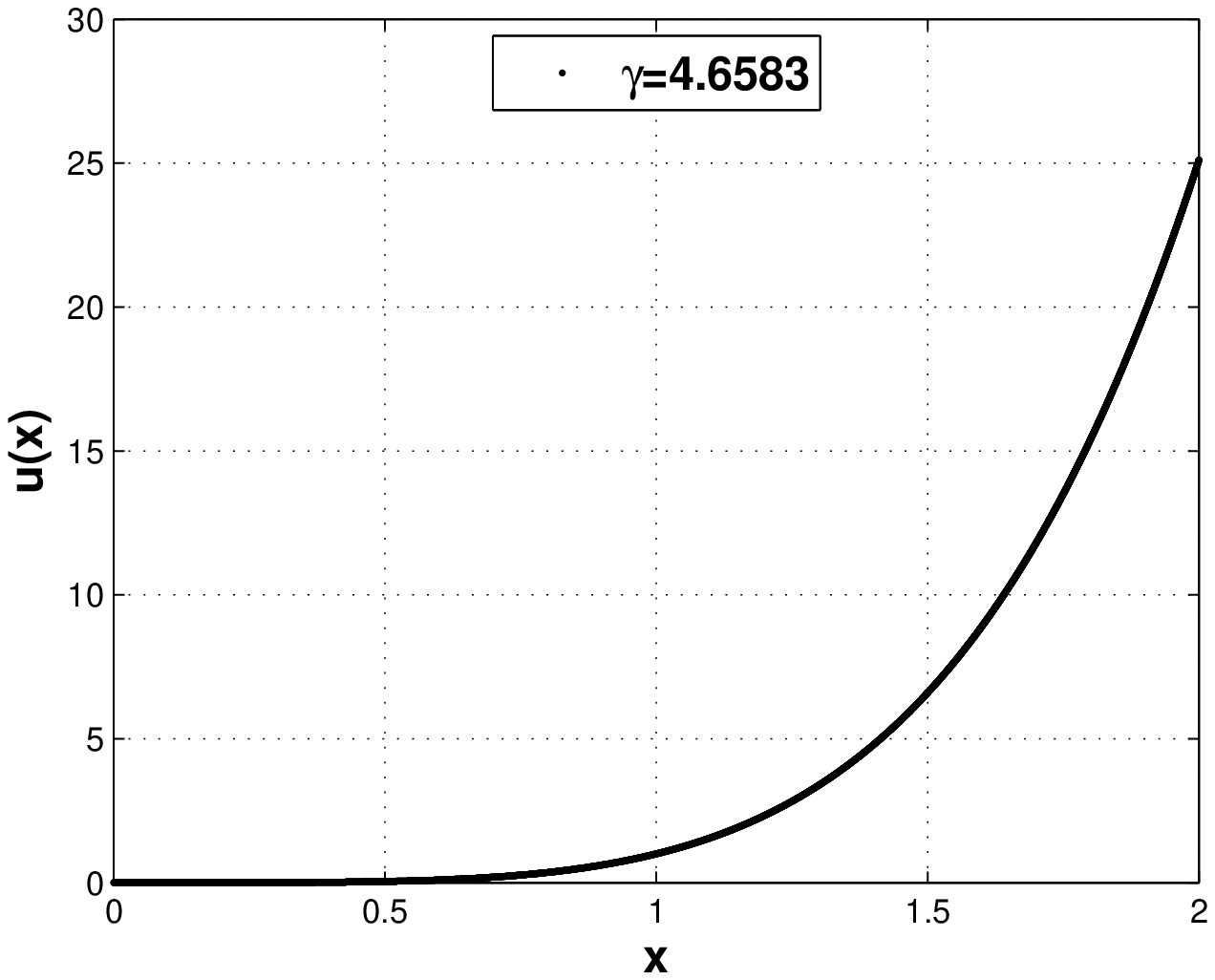}
		\includegraphics[width=6cm]{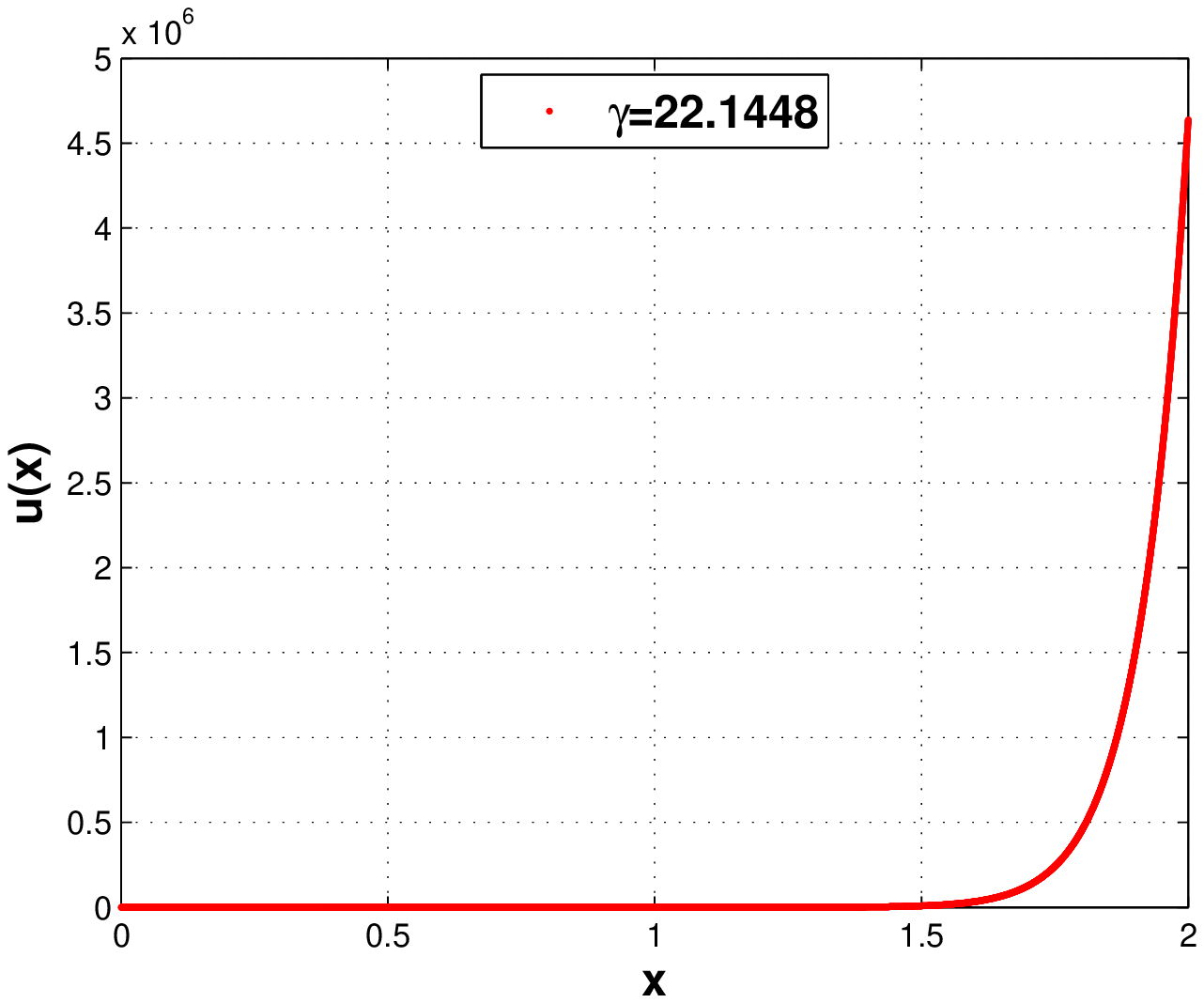}
		\caption{Solutions for equation \eqref{Example2b} with $\nu=1$.}	\label{Example2bGraph}
	\end{center}
\end{figure}

\vspace*{-8pt}

Hence, the solution in the form of a series is not unique.
\end{example}
\vspace*{-5pt}

\begin{example} \label{ex5} (no solutions with $d_i < 0$).\\ 
Let us consider another variation of equation \eqref{Example2a}, which also does not fit the constructed theory where all $d_i$ are positive,
\vskip -12pt 
\begin{equation}\label{Example2c}
   -0.1x^{4.9}D^{4.9}u(x)-6x^{3.75}D^{3.75}u(x)+x^{2.7}D^{2.7}u(x)+(x^{3.1}-0^2)u(x)=0.
\end{equation}
In this case the calculated restriction on $\nu^2_{\min}$ is negative, thus we have no restrictions on $\nu$. But the function $G(\gamma)-\nu^2$ has no zeros for \linebreak
$\gamma > n_{\max}-1 = 4$, because the largest root of equation \eqref{CondForC0} is \linebreak
$\gamma=2.9195$.
The increase in $\nu$ decreases  $\gamma$  -- the largest root of equation \eqref{CondForC0} -- until all solutions completely disappear. See Figure \ref{GyforEx5}.

\begin{figure}[H]
	\begin{center}
		\includegraphics[width=6cm]{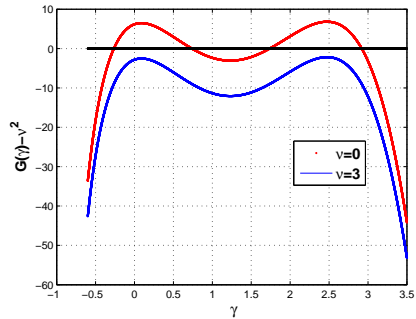}
		\caption{Function $G(\gamma)-\nu^2$ for $\nu=0$ and $\nu=3$, equation \eqref{Example2c}.}\label{GyforEx5}
	\end{center}
\end{figure}
\end{example}

\vspace*{-6pt}

\begin{example} \label{ex6} (one solution with $d_i < 0$).\\
Let us consider a slight variation of equation \eqref{Example2c}, which also does not fit the constructed theory but nevertheless produces one solution.
\begin{equation}\label{Example2d}
0.1x^{4.9}D^{4.9}u(x)-6x^{3.75}D^{3.75}u(x)+x^{2.7}D^{2.7}u(x)+(x^{3.1}-3^2)u(x)=0.
\end{equation}
In this case the calculated minimum $\nu^2_{\min}$ is also negative but there exists the large root $\gamma$ for equation \eqref{CondForC0}: $\gamma=38.8813 > n_{\max} - 1 = 4$.  Figure \ref{GyforEx6} demonstrates the solution for equation \eqref{CondForC0} with parameters generated by equation \eqref{Example2d}.  Figure \ref{Example6Graph} shows the solution of equation \eqref{Example2d}.

\begin{figure}[H]
	\begin{center}
		\includegraphics[width=6cm]{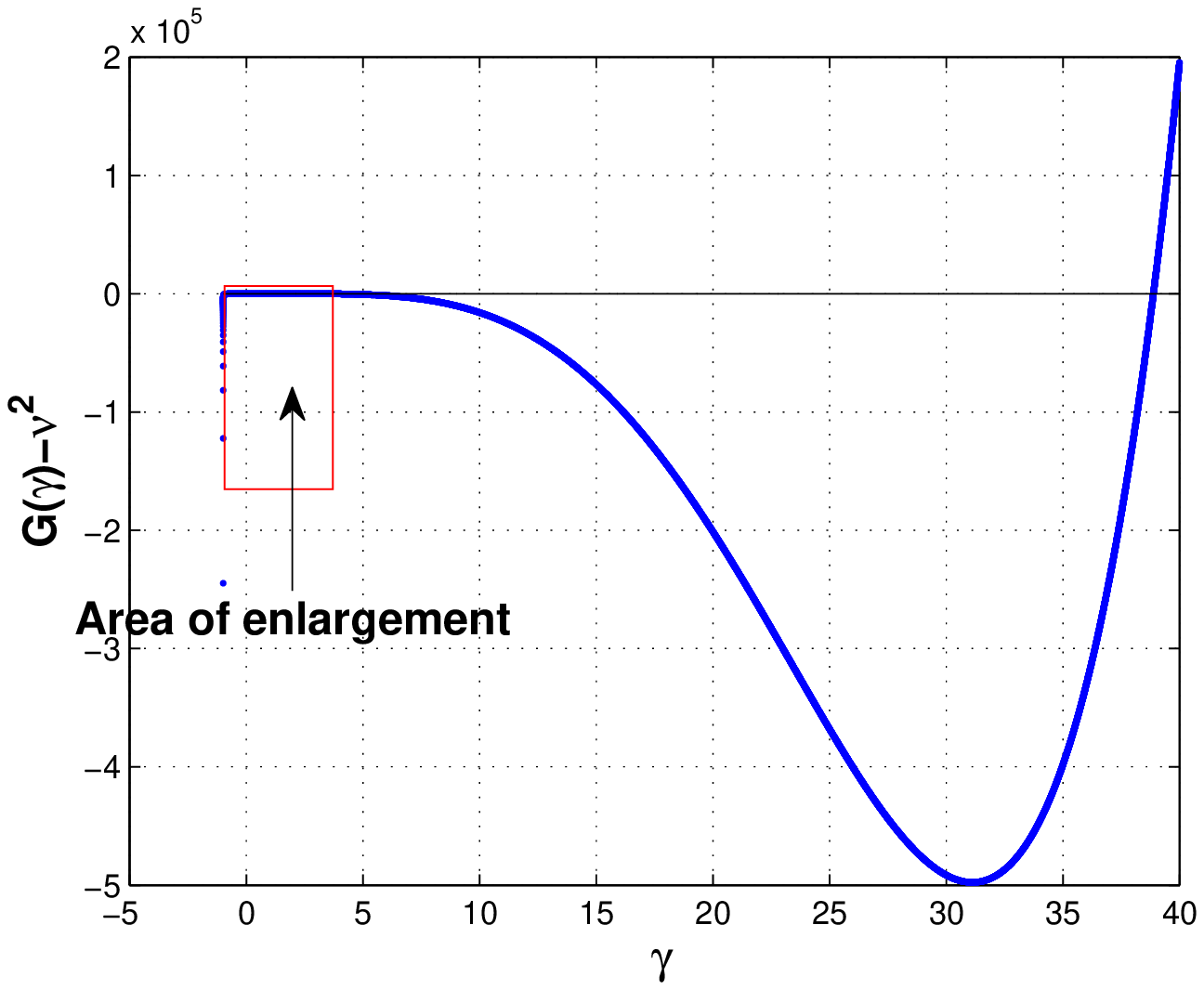}
		\includegraphics[width=6cm]{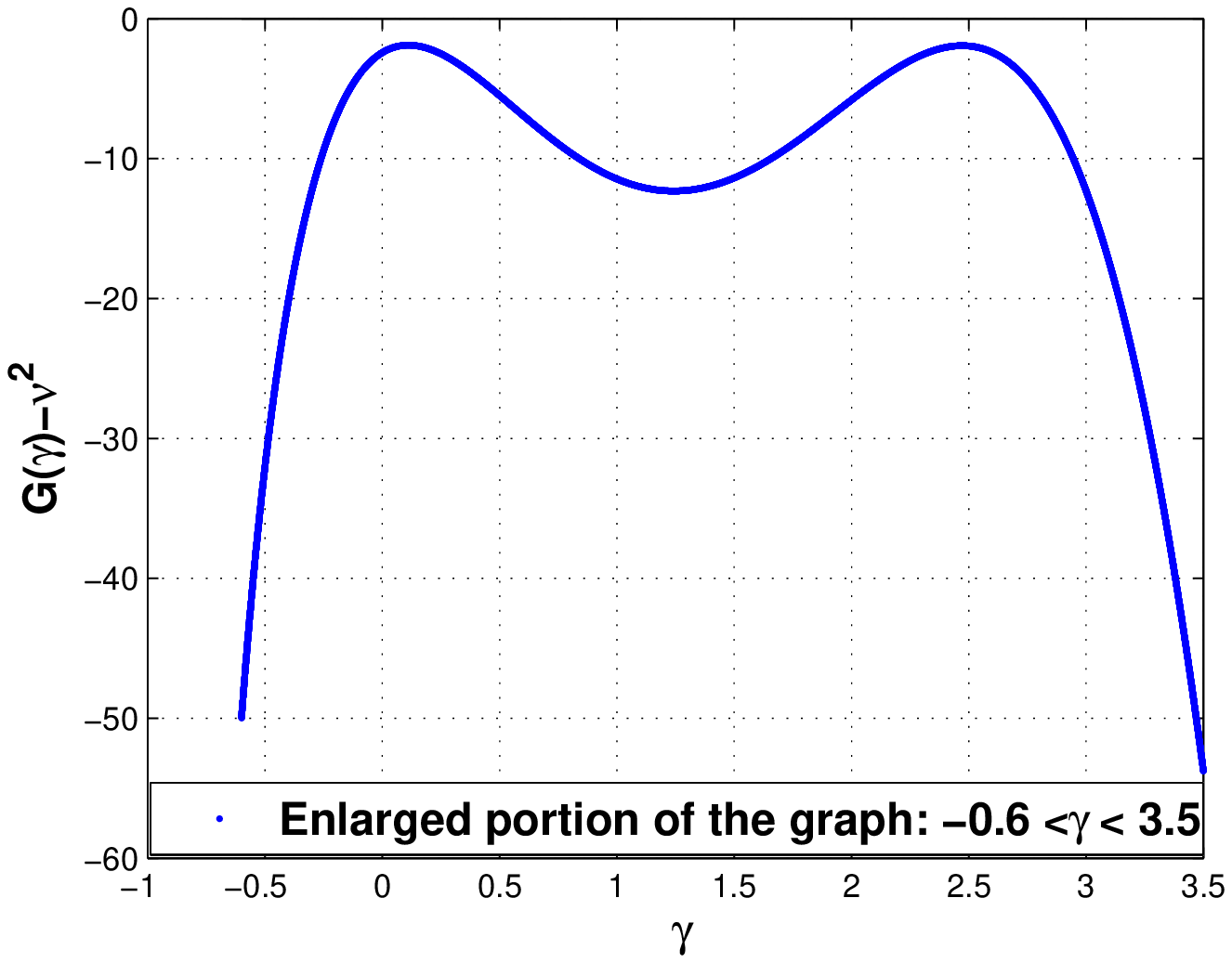}
		\caption{Function $G(\gamma)-3^2$ for equation \eqref{Example2d} for two intervals of $\gamma$.}\label{GyforEx6}
	\end{center}
\end{figure}

\vspace*{-10pt} 

\begin{figure}[H]
	\begin{center}
		\includegraphics[width=6cm]{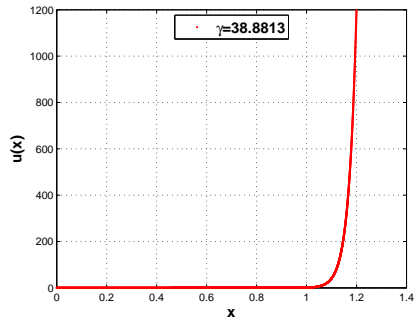}
		\caption{Solution $u(x)$ for equation \eqref{Example2d}.}\label{Example6Graph}
	\end{center}
\end{figure}

\vspace*{-3pt}

In this case one solution exists for any value of $\nu \in \mathbb{R}$.
\end{example}

\vspace*{-3pt}

\begin{example} \label{ex7} (multiple solutions with large integer derivative).\\
In equation \eqref{Example7} $\alpha_{\max}=6$, $n_{\max}=3$.
Since in this case $\alpha_{\max}-1 > n_{\max}$, then, in accordance with Theorem \ref{Thrm1},
the uniqueness condition fails.
\begin{equation}\label{Example7}
x^{6}u^{(6)}(x)+0.02x^{2.7}D^{2.7}u(x)+0.1x^{1.2}D^{1.2}u(x)+(x^{2}-1.5^2)u(x)=0.
\end{equation}

Function $G(\gamma)-\nu^2$ for which $\nu_{\min}^2=1.5671$ has for $\gamma > n_{\max}-1=2$ three acceptable
solutions
$\gamma_1=3.1395, \gamma_2=3.9392, \gamma_3=5.0056$ as we can see in Figure \ref{GyforEx7}.

\begin{figure}[H]
                \begin{center}
                         \includegraphics[width=6cm]{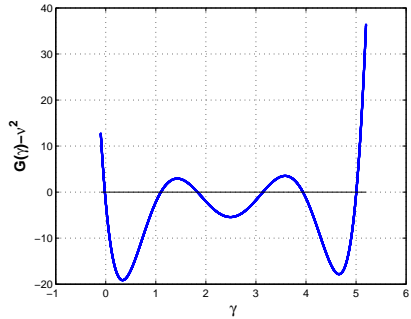}
                       \caption{Function $G(\gamma)-\nu^2$ for equation \eqref{Example7}.}              \label{GyforEx7}
                \end{center}
\end{figure}

These three solutions for equation \eqref{Example7} are presented in Figures \ref{Example7Graph},\ref{Example7Graph-2}. The found values
$\gamma< n_{\rm max}-1=2$ do not produce solutions to the differential equation.

\begin{figure}[H]
     \begin{center}
            \includegraphics[width=6cm]{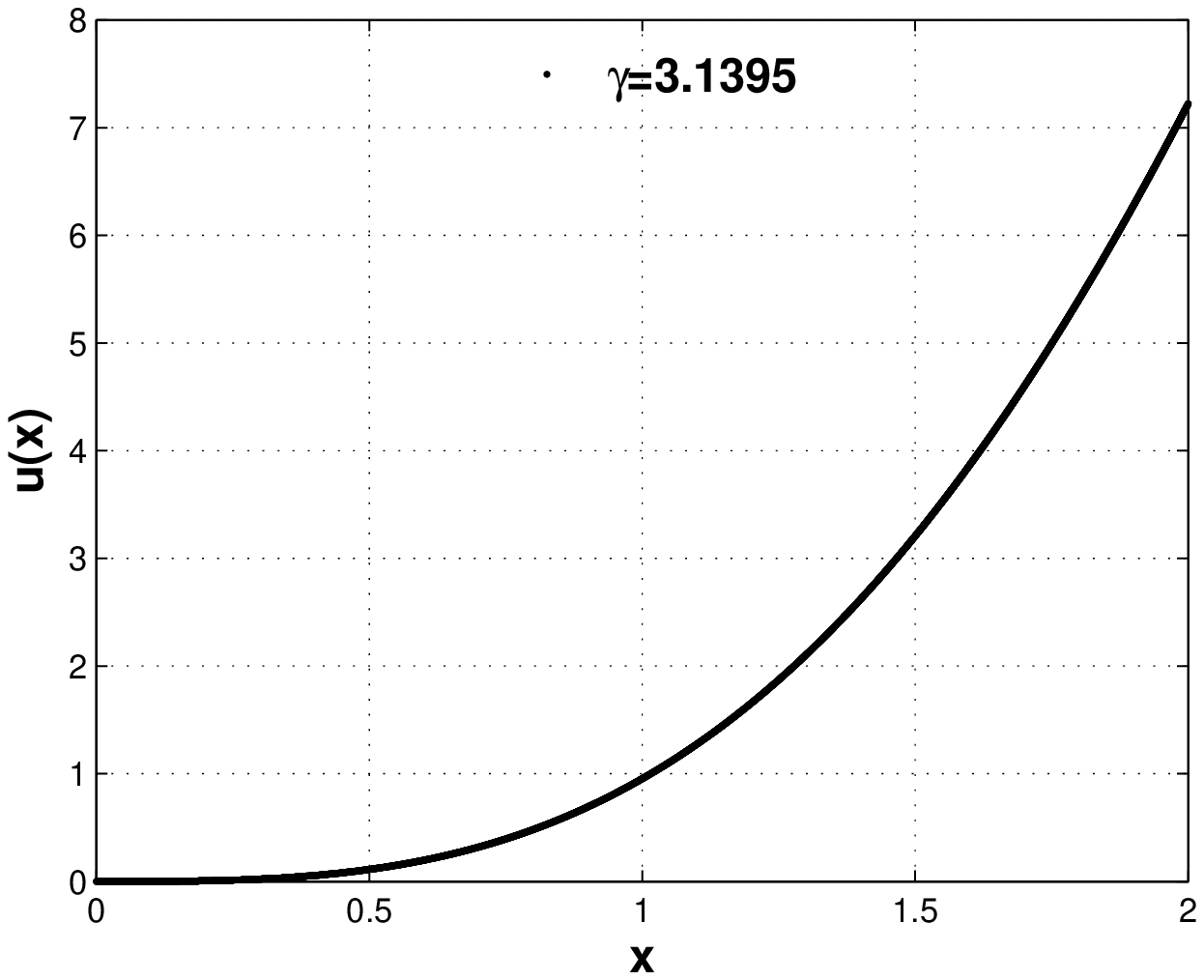}
            \includegraphics[width=6cm]{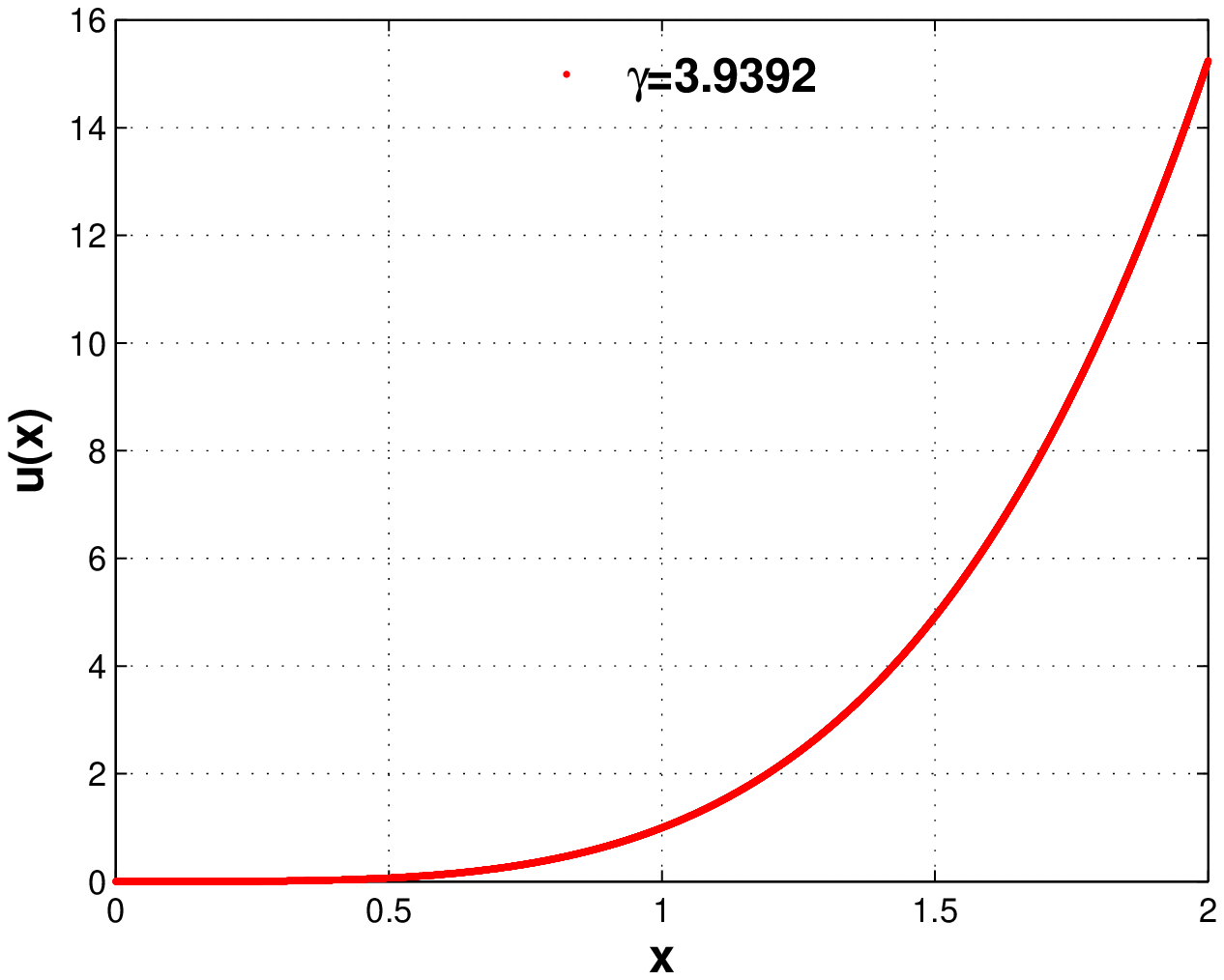}
          \caption{Solutions one and two for equation \eqref{Example7}.}\label{Example7Graph}
     \end{center}
\end{figure}

\vspace*{-8pt}

\begin{figure}
	\begin{center}
          \includegraphics[width=6cm]{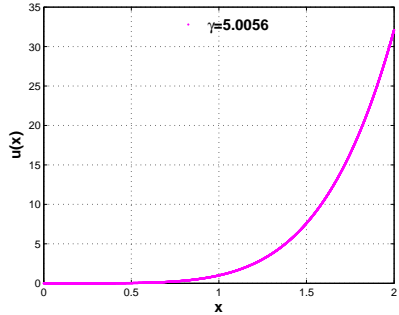}
          \caption{Third solution for equation \eqref{Example7}.}\label{Example7Graph-2}
    \end{center}
\end{figure}
\end{example}

\vspace*{-3pt}%

\begin{example} \label{ex8} (multiplicity of $\gamma$ is 2)
Are there any equations such that the generated characteristic equation has multiplicity greater than 1?  The examples below show that such equations definitely exist.
Let us consider equation
\begin{equation}\label{Example8Sect3}
-2x^{1.9}D^{1.9}u(x)+10x^{0.9}D^{0.9}u(x)-5x^{0.7}D^{0.7}u(x)+(x-6.336632736437)=0.
\end{equation}
In this case $\nu=2.517266917996$.  The graph of $G(\gamma)=\displaystyle\sum_{i=1}^{m}\frac{d_i\Gamma(1+\gamma)}{\Gamma(1+\gamma-\alpha_i)}$  (see \eqref{CondForC0}) is given in Figure \ref{GyforEx8Sect3}.

\begin{figure}[H]
	\begin{center}
		\includegraphics[width=6cm]{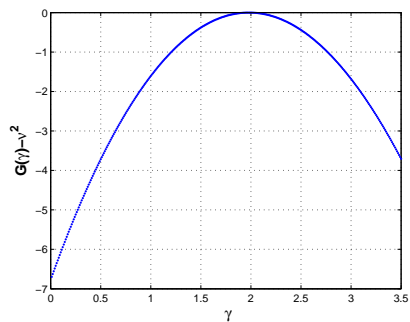}
		\caption{Function $G(\gamma)-\nu^2$ for equation \eqref{Example8Sect3}.}	\label{GyforEx8Sect3}
	\end{center}
\end{figure}

At $\gamma=1.979$ the derivative of the graph is zero and the multiplicity of the root is 2.  In this case we obtain two linearly independent solutions as shown in Figure \ref{Ex8solnsSect3}.

\begin{figure}[H]
	\begin{center}
		\includegraphics[width=6cm]{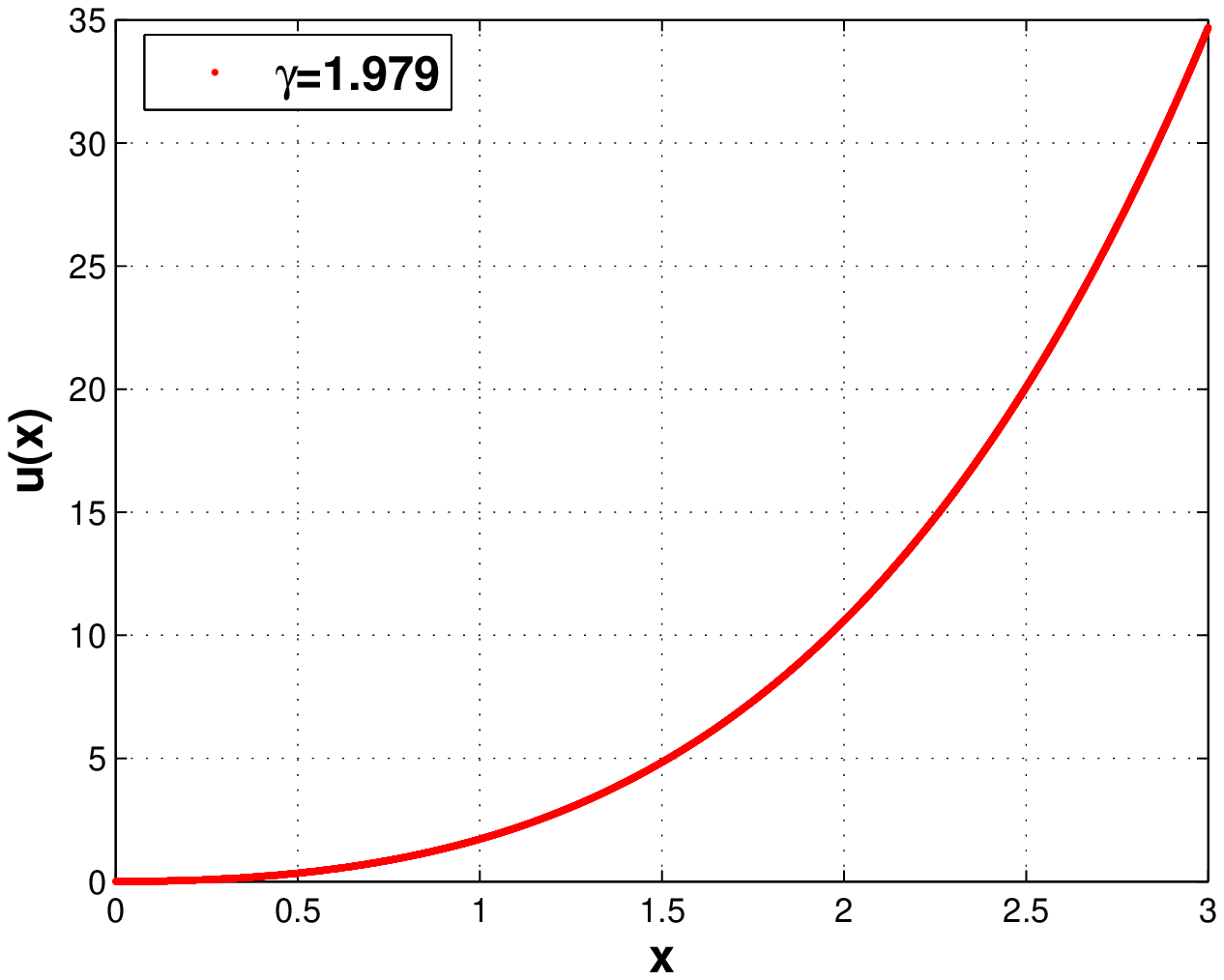}
		\includegraphics[width=6cm]{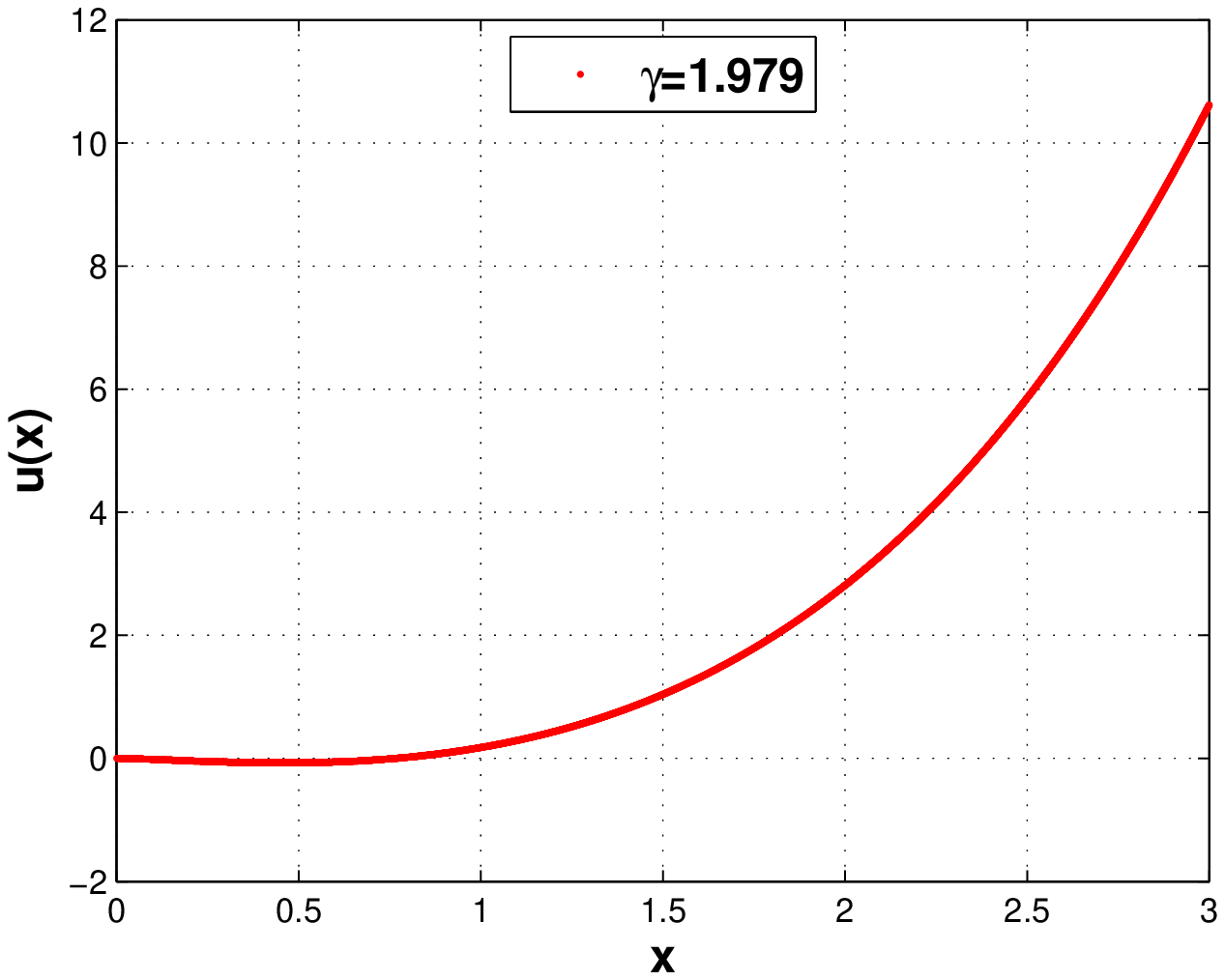}
		\caption{Solutions for equation \eqref{Example8Sect3}.}\label{Ex8solnsSect3}
	\end{center}
\end{figure}
\end{example}
\vspace*{-8pt}

\section*{Conclusions}

\begin{itemize}
\item We introduce the generalized fractional Bessel equation, which may include integer derivatives and covers the fractional and classical Bessel equations as particular cases.
\item We construct the solution in the form of fractional power series for the generalized fractional Bessel equation.
\item The existence and uniqueness theorems for the series solution have been proved.
\item Solutions to the equations when multiplicity of the root of the characteristic equation is greater than one are identified.
\item The uniqueness theorem for the initial value problem in the space of continuously differentiable functions has been proved.
\item Several numerical examples support the constructed theory and provide the cases when the infringement of the conditions of the theorem leads to non-existence or non-uniqueness of the solution and, thus, these counterexamples justify the necessity of our requirements.
\end{itemize}

\vspace*{-2pt}

\section*{\small Acknowledgements}

The authors would like to thank Professor L.~Boyadjiev for drawing our attention to the fractional Bessel equation and Professor V.~Kiryakova for the valuable remarks, which helped to improve the paper.


\end{document}